\documentclass[11pt]{article}
\usepackage{amsfonts,amssymb}
\usepackage{diagrams}

\usepackage[latin1]{inputenc}  
\usepackage[T1]{fontenc}

\newcounter{paragrafsubsub}[subsubsection]
\renewcommand{\theparagrafsubsub}{%
\thesubsubsection.\roman{paragrafsubsub}}
\newcommand{\paragrafsubsub}{%
\refstepcounter{paragrafsubsub}
{\bf \theparagrafsubsub}\hspace{0.2em}--- }

\newcounter{paragrafsub}[subsection]
\renewcommand{\theparagrafsub}{\thesubsection.\arabic{paragrafsub}}
\newcommand{\paragrafsub}{%
\refstepcounter{paragrafsub}
{\bf \theparagrafsub}\hspace{0.2em}--- }

\newcounter{paragraf}[section]
\renewcommand{\theparagraf}{\thesection.\arabic{paragraf}}
\newcommand{\paragraf}{%
\refstepcounter{paragraf}
{\bf \theparagraf}\hspace{0.2em}--- }

\newcommand\paragraphe{%
\par \indent
\ifcase\value{subsection} %
\paragraf
\else
\ifcase\value{subsubsection}\paragrafsub %
\else\paragrafsubsub
\fi\fi
}

\def\longto{\longrightarrow}

\def\cone{{\mathcal C}}
\def\tc{{\mathcal{TC}}}
\def\ac{{\mathcal{AC}}}
\def\sac{{\mathcal{SAC}}}
\def\lr{{\mathcal{LR}}}
\def\Chi{X}\def\Y{Y}
\def\Det{{\mathcal Det}}
\def\Tau{{\mathcal T}}
\def\Hcal{{\mathcal H}}

\def\PP{{\mathbb P}}\def\GG{{\mathbb G}}
\def\RR{{\mathbb R}}\def\QQ{{\mathbb Q}}\def\ZZ{{\mathbb Z}}

\def\kk{{\mathbb K}}

\def\lg{{\mathfrak g}}

\def\lp{{\mathfrak p}}\def\lb{{\mathfrak b}}\def\lu{{\mathfrak u}}

\def\hL{{\hat L}}
\def\hlg{{\hat{\lg}}}\def\hlb{{\hat{\lb}}}\def\hlu{{\hat{\lu}}}\def\hlp{{\hat{\lp}}}

\def\Lim{{p_\lambda}}
\def\kbprod{\odot_0}

\def\Orb{{\mathcal O}}
\def\H{{\rm H}}

\def\Li{{\mathcal{L}}}
\def\Mi{{\mathcal{M}}}
\def\O{{\mathcal O}}
\def\Oun{{\mathcal O}(1)}

\def\CGLX{\ac^G_\Lambda(X)}

\def\quot{/\hspace{-.5ex}/}
\def\Om{\Omega}

\def\GPY{G\times_PY}

\def\Face{{\mathcal F}}

\def\St{\rm St}
\def\Pic{\rm Pic}

\def\codim{{\rm codim}}

\def\hG{{\hat{G}}}\def\hg{{\hat{g}}}
\def\hB{{\hat{B}}}\def\hU{{\hat{U}}}
\def\hT{{\hat{T}}}
\def\hW{{\hat{W}}}\def\hw{{\hat{w}}}
\def\hnu{{\hat\nu}}\def\htheta{{\hat\theta}}\def\hchi{{\hat\chi}}
\def\hP{{\hat P}}\def\hQ{{\hat Q}}
\def\hR{{\hat R}}
\def\hrho{{\hat{\rho}}}

\def\P{{\mathcal P}}\def\Q{{\mathcal Q}}

\def\Im{{\rm Im}}

\newtheorem{lemma}{Lemma}
\newtheorem{prop}{Proposition}
\newtheorem{theo}{Theorem}
\newtheorem{coro}{Corollary}

\newtheorem{theoi}{Theorem}

\newenvironment{proof}{{\noindent\bf Proof.}}{\hfill $\square$}
\newenvironment{defin}{{\noindent\bf Definition.}}{\\}
\newenvironment{remark}{{\noindent\bf Remark.}}{}
\newenvironment{nota}{{\noindent\bf Notation.}}{\\}

\begin{document}
\title{Geometric Invariant Theory and\\
Generalized Eigenvalue Problem}
\author{N. Ressayre\footnote{Universit{\'e} Montpellier II - 
CC 51-Place Eug{\`e}ne Bataillon -
34095 Montpellier Cedex 5 -
France - {\tt ressayre@math.univ-montp2.fr}}}

\maketitle

\begin{abstract}
Let $G$ be a connected reductive subgroup of a complex connected reductive group $\hat{G}$.
Fix maximal tori and Borel subgroups of $G$ and $\hG$.
Consider the cone $\lr(G,\hG)$ generated by the pairs $(\nu,\hnu)$ of dominant characters 
such that $V_\nu$  is a submodule of $V_{\hnu}$ (with usual notation).
Here we give a minimal set of inequalities describing $\lr(G,\hG)$ as a part of the dominant 
chamber. 

In way, we obtain results about the faces of the Dolgachev-Hu's $G$-ample cone and variations 
of this cone. 
\end{abstract}

\section{Introduction}

In this article, we are mainly interested in the faces of the $G$-ample cone
as defined by Dolgachev-Hu in \cite{DH}.
In this introduction we start by the end: we firstly explain the applications
to the generalized  Horn's problem.\\ 

Let $G$ be a connected reductive subgroup of a connected reductive group $\hat{G}$ both defined over an algebraic
closed field $\kk$ of characteristic zero. 
We consider the following question:

\begin{center}
  What irreducible representations of $G$ appear in a given irreducible representation of $\hG$ ?
\end{center}

Once maximal tori ($T\subset\hT)$ and Borel subgroups ($B\supset T$ and $\hB\supset \hT$) fixed, 
the question  is to understand the set ${\rm LR}(G,\hG)$ of pairs $(\nu,\hnu)$ of dominant 
character of $T\times\hT$ such that the $G$-module associated to $\nu$ can be $G$-equivariantly 
embedded in the $\hG$-module associated to $\hnu$. 
By a result of M.~Brion and F.~Knop (see \cite{elash}), ${\rm LR}(G,\hG)$ is a 
finitely generated submonoid of the character group of $T\times\hT$. 
Our purpose is to study the linear inequalities satisfied by this monoid. 
More precisely, we consider the convex cone $\lr(G,\hG)$ generated by 
 ${\rm LR}(G,\hG)$:
this cone is characterized by finitely many inequalities.
It occurs that the interior of $\lr(G,\hG)$ is non empty if and only if no non trivial 
connected normal subgroup
of $G$ is normal in $\hG$: for simplicity, we assume, from now on that $\lr(G,\hG)$ has non empty interior.
Our first result is a list of inequalities  which characterize $\lr(G,\hG)$ as a part
of the dominant cone. Such a list was already obtained in \cite{sjamaar}. 
Whereas the Berenstein-Sjamaar's list is redundant, our list is proved to be irredundant.\\

To make our statements more precise, we  introduce notation.
Consider the natural paring  $\langle\cdot\,,\,\cdot\rangle$ between one parameter subgroups and 
characters of tori $T$ or $\hT$.
Let $W$ (resp. $\hW$) denote the Weyl group of $T$ (resp. $\hT$).
If $\lambda$ is a one parameter subgroup of $T$ (or so of $\hT$), 
we denote by $W_\lambda$ (resp. $\hW_\lambda$) the stabilizer of $\lambda$ for the natural action
of the Weyl group on the set of one parameter subgroups.
For $w\in W/W_\lambda$ and $\hat{w}\in \hW/\hW_\lambda$, we consider the following linear form on the 
character group $\Chi(T\times\hT)$ of $T\times \hat T$:
$$
\varphi_{\lambda,w,\hw}\,:\,
(\nu,\hat\nu)\mapsto\langle \hat w\lambda,\hat\nu\rangle+\langle w\lambda,\nu\rangle.
$$
In fact, all the inequalities in Theorem~\ref{thi:BeSj} below have the following form
$\varphi_{\lambda,w,\hw}\geq 0$. We need some notation to explain which triples $(\lambda,w,\hw)$ appear.

Let $P(\lambda)$ denote  the usual parabolic subgroup of $G$ associated to $\lambda$ 
(see Section~\ref{sec:aops}).
The cohomology group ${\rm H}^*(G/P(\lambda),\ZZ)$ is freely generated by the Schubert classes
$[\Lambda_w]$ parametrized by the  elements $w\in W/W_\lambda$.
We will consider $\hG/\hP(\lambda)$, $[\Lambda_\hw]$ as above but with $\hG$ in 
place of $G$.
Consider also the canonical $G$-equivariant immersion 
$\iota\,:\,G/P(\lambda)\longto\hG/\hP(\lambda)$; and the corresponding morphism 
in cohomology $\iota^*$.

Let $\lg$ and $\hlg$ denote the Lie algebras of $G$ and $\hG$.
Let $\rho$, $\hat\rho$ and $\hat\rho^\lambda$ be  the half-sum of positive roots of $G$, $\hat G$ and 
the centralizer $\hG^\lambda$ of $\lambda$ in $\hat G$.
Let $\lambda_1,\,\cdots,\lambda_n$ be the set of indivisible dominant one parameter subgroups of
$T$ orthogonal to an hyperplane generated by weights of $T$ in $\hlg/\lg$.

We can now state one of our main results:

\begin{theoi}\label{thi:BeSj}
    We assume that no ideal of $\lg$ is an ideal of $\hlg$. 
Then, $\lr(G,\hG)$ has non empty interior  in $\Chi(T\times\hT)\otimes \QQ$.

A dominant weight $(\nu,\hnu)$ belongs to  $\lr(G,\hG)$ if and only if
for all $i=1,\cdots, n$ 
and for all pair of Schubert classes $([\Lambda_w],\,[\Lambda_\hw])$ of $G/P(\lambda_i)$ and 
$\hG/\hP(\lambda_i)$ associated to a pair $(w,\hw)\in W/W_{\lambda_i}\times \hW_{\lambda_i}$ such that 
 
\begin{enumerate}
\item \label{cond:BeSj1}
$\iota^*([\Lambda_\hw]).[\Lambda_w]=[{\rm pt}]\in H^*(G/P(\lambda),\ZZ)$, and 
\item  \label{cond:BeSj2}
$\langle w\lambda_i,\rho\rangle+\langle \hw\hat\lambda_i,\hat\rho\rangle=
\langle \lambda_i,\rho\rangle+\langle\lambda_i,2\hat\rho^{\lambda_i}-\hat\rho\rangle $,
\end{enumerate}
we have
\begin{eqnarray}
  \label{eq:ineg}
\langle w\lambda_i, \nu \rangle+\langle\hw\lambda_i,\hnu\rangle  \geq 0.
\end{eqnarray}
\end{theoi}

In \cite{sjamaar}, Berenstein and Sjamaar showed that $(\nu,\hnu)$ belongs to  $\lr(G,\hG)$ 
if and only if
$\langle w\lambda_i, \nu \rangle+\langle\hw\lambda_i,\hnu\rangle \geq 0$
for all $i=1,\cdots, n$ 
and for all pair of Schubert classes $([\Lambda_w],\,[\Lambda_\hw])$  such that
$\iota^*([\Lambda_\hw]).[\Lambda_w]=d.[{\rm pt}]$ for some positive 
integer $d$. 
In the case when $G$ is diagonally embedded in $G^s$, Kapovich-Leeb-Millson already proved that one may assume
that $d=1$. Again in the case when $G$ is diagonally embedded in $G^s$, Belkale and Kumar obtained 
the same inequalities as in  Theorem~\ref{thi:BeSj} in \cite{BK}.\\

In some sense, our second main result asserts that Theorem~\ref{thi:BeSj}  is optimal:

\begin{theoi}
\label{thi:mini} 
In Theorem~\ref{thi:BeSj}, Inequalities~(\ref{eq:ineg}) are pairwise distinct and no one can be omitted. 
\end{theoi}

This result was known is some particular cases.
Indeed, Knutson, Tao and Woodward shown in \cite{KTW} the case when $G={\rm SL}_n$ is diagonally embedded 
in ${\rm SL}_n\times{\rm SL}_n$ using combinatorial tools. 
Using the interpretation of the Littlewood-Richardson coefficients  
as structure coefficient of the cohomology ring of the Grassmann varieties,
Belkale made a geometric proof of Knutson-Tao-Woodward's result (see~\cite{belkale:geomHorn}).
Using explicit calculation with the help of a computer, Kapovich, Kumar
and Millson proves the case when $G={\rm SO}(8)$ is diagonally embedded in $G\times G$ in \cite{KKM}.
Our proof is different and uses Geometric Invariant Theory.\\

Let $([\Lambda_w],\,[\Lambda_\hw])$ be a  pair of Schubert classes  in $G/P(\lambda_i)$ such that
$\iota^*([\Lambda_\hw]).[\Lambda_w]=d.[{\rm pt}]\in H^*(G/P(\lambda),\ZZ)$ for some positive 
integer $d$.  By \cite{sjamaar}, the set of  $(\nu,\hnu)\in \lr(G,\hG)$ 
such that 
$\langle w\lambda_i, \nu \rangle+\langle\hw\lambda_i,\hnu\rangle=0$
is a face of  $\lr(G,\hG)$. 
Theorem~\ref{thi:mini} shows that if $(w,\hw)$ does not satisfy Conditions~\ref{cond:BeSj1} and \ref{cond:BeSj2} of
Theorem~\ref{thi:BeSj} then the codimension of this face is greater than one. 
We can improve this statement:

\begin{theoi}
\label{thi:redundantineq}
Let us fix a $\lambda_i$. Let  $([\Lambda_w],\,[\Lambda_\hw])$ be a  pair of Schubert varieties  in $G/P(\lambda_i)$ such that
$\iota^*([\Lambda_\hw]).[\Lambda_w]=d.[{\rm pt}]\in H^*(G/P(\lambda),\ZZ)$ for some positive integer $d$.
We assume that there exists a pair of dominant weights $(\nu,\hnu)\in\lr(G,\hG)$ such that
$\langle w\lambda_i, \nu \rangle+\langle\hw\lambda_i,\hnu\rangle=0$ and such that $\nu$ is strictely dominant.

Then, $d=1$ and $\langle w\lambda,\rho\rangle+\langle \hw\hat\lambda,\hat\rho\rangle=
\langle \lambda,\rho\rangle+\langle\lambda,2\hat\rho^\lambda-\hat\rho\rangle $.

In particular, the set of $(\nu,\hnu)\in\lr(G,\hG)$ such that
$\langle w\lambda_i, \nu \rangle+\langle\hw\lambda_i,\hnu\rangle=0$
is a face of codimension one.
\end{theoi}

Theorems~\ref{thi:BeSj} and \ref{thi:mini} can be thought as a description of the faces of codimension one of 
$\lr(G,\hG)$ which intersect the interior of the dominant chamber. 
In Theorem~\ref{th:ppal} below, we study the smaller faces of $\lr(G,\hG)$.
To avoid some notation in the introduction, we only state our results about these 
faces in the case when $G$ is diagonally embedded in $G^s$ for an integer $s\geq 2$; 
that is, in the case of decomposition of the tensor product (see Section~\ref{sec:GBGB} for a 
more general result).

If $\alpha$ is a simple root of $G$, $\omega_{\alpha^\vee}$ denote the corresponding fundamental 
weight. If $I$ is a set of simple roots, $P(I)$ denote the standard parabolic subgroup associated to
$I$ and $W_I$ its Weyl group.
In  \cite{BK}, Belkale and Kumar defined a new product $\kbprod$ on the cohomology groups 
${\rm H}^*(G/P,\ZZ)$.

\begin{theoi}
\label{thi:petiteface}
We assume that $G$ is semisimple diagonally embedded in $\hG=G^s$ for an integer $s\geq 2$.

\begin{enumerate}
\item  Let $I$ be a set of $d$ simple roots
and $(w_1,\cdots,w_{s+1})\in (W/W_I)^{s+1}$  such that
$[\Lambda_{w_1}]\kbprod\cdots\kbprod[\Lambda_{w_{s+1}}]=[\Lambda_e]\in H^*(G/P(I),\ZZ)$.
Then, the set of $(\nu_1,\cdots,\nu_{s+1})\in \lr(G,G^s)$ such that 
$$
\forall \alpha\in I\hspace{3ex}\sum_i\langle\omega_{\alpha^\vee},w_i^{-1}\nu_i \rangle =0,
$$
is a face of codimension $d$ of $\lr(G,G^s)$.
\item Any face of $\lr(G,G^s)$ intersecting the interior of the dominant chamber of $G^{s+1}$ 
is obtained in this way.
\end{enumerate}
\end{theoi}

If there is a lot of literature on the faces of codimension one, I do not know any other result
about smaller faces even in the case of Horn's conjecture that is for 
${\rm SL}_n\subset{\rm SL}_n\times{\rm SL}_n$.
Moreover, Theorem~\ref{thi:petiteface} gives an application of Belkale-Kumar's product $\kbprod$
for any $G/P$ whereas in \cite{BK} only the case when $P$ is maximal is used.\\

Let us explain the starting point of the proofs of Theorems~\ref{thi:BeSj} to \ref{thi:petiteface}.
Consider the variety $X=G/B\times \hG/\hB$ endowed with the diagonal action of $G$.
To any character $(\nu,\hnu)$ of $T\times \hT$, one associates a $G$-linearized line bundle
$\Li_{(\nu,\hnu)}$ on $X$ such that ${\rm H}^0(X,\Li_{(\nu,\hnu)})=V_\hnu\otimes V_\nu$. 
Hence, $(\nu,\hnu)$ belongs to $\lr(G,\hG)$ if and only if a positive power of $\Li_{(\nu,\,\hnu)}$ 
admit a non zero $G$-invariant section. If $\Li_{(\nu,\,\hnu)}$ is ample this is equivalent to saying 
that $X$ admit semistable points for $\Li_{(\nu,\,\hnu)}$. So, one can use classical results of Geometric
Invariant Theory as Hilbert-Mumford's Theorem and Luna's Slice Etale Theorem.
This method was already used in \cite{Kly,sjamaar,BK}.\\

We obtain results in the following more general context. 
Consider a connected reductive group $G$ acting on a normal projective variety $X$.
To any $G$-linearized line bundle $\Li$ on $X$ we associate the following open subset 
$X^{\rm ss}(\Li)$  of $X$:
$$
X^{\rm ss}(\Li)=\left \{
x\in X \mbox{ : }
\exists n>0\mbox{ and }
\sigma\in{\rm H}^0(X,\Li^{\otimes n})^G \mbox{ such that }\sigma(x)\neq 0
\right \}.
$$
The points of $X^{\rm ss}(\Li)$ are said to be {\it semistable} for $\Li$.
Note that if $\Li$ is not ample, this notion of semistability is not the standard one. 
In particular, the quotient $\pi_\Li\::\:X^{\rm ss}(\Li)\longto X^{\rm ss}(\Li)\quot G$ is a 
good quotient, if $\Li$ is ample; but not in general.
In this context, we ask for:
\begin{center}
  What are the $\Li$'s with non empty set $X^{\rm ss}(\Li)$ ?
\end{center}
Let us fix a freely finitely generated subgroup $\Lambda$ of the group $\Pic^G(X)$ of $G$-linearized
line bundles on $X$. 
Let $\Lambda_\QQ$ denote  the $\QQ$-vector space containing $\Lambda$ as a lattice.
Consider the convex cones $\tc_\Lambda^G(X)$ (resp. 
$\ac_\Lambda^G(X)$) generated in $\Lambda_\QQ$ by the $\Li$'s (resp. the ample $\Li$'s) in $\Lambda$ 
which have non zero $G$-invariant sections. By~\cite{DH} (see also \cite{GeomDedic}), 
$\ac^G_\Lambda(X)$ is a closed convex rational polyhedral cone in the dominant cone of $\Lambda_\QQ$.
We are interested in the faces of $\ac^G_\Lambda(X)$ and $\tc_\Lambda^G(X)$.\\

\bigskip
We are now going to define the notion of well covering pair which will play a central role
in this work.
Let $\lambda$ be a one parameter subgroup of $G$ and $C$ be an irreducible 
component of its fix points.
Consider $C^+=\{x\in X\ |\ \lim_{t\to 0}\lambda(t)x\in C\}$ and
the natural $G$-equivariant map $\eta\::\:G\times_{P(\lambda)}C^+\longto X$. 
The pair $(C,\lambda)$ is said to be {\it well covering} if $\eta$ induces an isomorphism over
an open subset of $X$ intersecting $C$. The pair is said to be {\it dominant} if $\eta$ is
dominant.

For $\Li\in\Pic^G(X)$, we denote by $\mu^\Li(C,\lambda)$ the integer giving the action of 
$\lambda$ on the restriction of $\Li$ on $C$. 
We now state a first description of $\ac^G_\Lambda(X)$:

\begin{theoi}
\label{thi:desc1}
An ample $G$-linearized line bundle $\Li\in\Lambda$ belongs to $\ac^G_\Lambda(X)$ 
if and only if  for all well covering pair
$(C,\lambda)$ with a dominant one parameter subgroup $\lambda$ of $T$ we have
$\mu^\Li(C,\lambda)\leq 0$.
\end{theoi}

Let $(C,\lambda)$ be a well covering pair. 
Theorem~\ref{thi:desc1} shows that the set of $\Li\in\CGLX$ such that 
$\mu^\Li(C,\lambda)= 0$ is a face of $\CGLX$.
We obtained several relations between the faces of $\CGLX$ and 
pairs $(C,\lambda)$.
Firstly, for $\Li$ in the boundary of $\CGLX$, we  describes the quotient
variety $X^{\rm ss}(\Li)\quot G$ in terms of covering pair:

\begin{theoi}
\label{thi:CsurGlambda}
We assume that $X=G/B\times Y$, for a normal projective $G$-variety $Y$.
Let $(C,\lambda)$ be a dominant pair and $\Li$ be an ample $G$-linearized line bundle on $X$ such 
that $\mu^\Li(C,\lambda)=0$.
Consider the action of the centralizer $G^\lambda$ of $\lambda$ in $G$ on $C$ and 
the associated set of semistable points $C^{\rm ss}(\Li,G^\lambda)$. 
We assume that $X^{\rm ss}(\Li)$ is not empty.

Then,
\begin{enumerate}
\item $X^{\rm ss}(\Li)\quot G$ is isomorphic to $C^{\rm ss}(\Li,G^\lambda)\quot G^\lambda$; and,
\item $(C,\lambda)$ is well covering.
\end{enumerate}
\end{theoi}

Actually, Theorem~\ref{thi:CsurGlambda} is a key point in the proof of Theorem~\ref{thi:redundantineq}.\\

We already noticed that to any dominant pair, one can associate a face (eventually empty) of $\CGLX$.
Conversely, we have:

\begin{theoi}
\label{thi:face2cov}
Let $\Face$ be a face of $\CGLX$.

Then, there exists a well covering pair $(C,\lambda)$ such that 
 $\Face$ is the set of the $\Li$'s in $\CGLX$ such that $\mu^\Li(C,\lambda)=0$.
\end{theoi}

We now assume that $\Lambda\otimes \QQ=\Pic^G(X)\otimes\QQ$. Notice that if $\Pic(X)$ is not finitely
generated it can be replaced by the Neron-Severi group (see~\cite{DH}).
We denote $\tc^G_\Lambda(X)$ by $\tc^G(X)$.
We are interested in a kind of converse of Theorem~\ref{thi:face2cov}.
So, we fix a well covering pair $(C,\lambda)$. 

\begin{theoi}
\label{thi:rec}
With above notation, consider the  linear map $\rho$ induced by the restriction:
$$
\rho\ :\ \Pic^G(X)\otimes\QQ\longto\Pic^{G^\lambda}(C)\otimes\QQ.
$$
Then, the subspace of $\Pic^G(X)\otimes\QQ$ spanned by the $\Li\in\tc^G(X)$ such that
$\mu^\Li(C,\lambda)=0$ is the pullback by $\rho$ of the intersection of the image of $\rho$ 
and $\tc^{G^\lambda}(C)$.  
\end{theoi}

By induction, Theorem~\ref{thi:rec} will allow us to prove Theorem~\ref{thi:mini}.
Indeed, it will be used to prove that the set of $\Li\in\lr(G,\hG)$ which realize a given equality in 
Theorem~\ref{thi:BeSj} is of codimension one.\\

In this article, we are mainly interested in the cones $\lr(G,\hG)$ which are examples of 
cones $\sac_\Lambda^G(X)$. Other interesting examples will be studied elsewhere. 
For example, if $Y$ is any $G$-variety endowed with an $G$-linearized line bundle $\Li$ 
the moment polytope $P(Y,\Li)$ is an affine section of a cone $\sac_\Lambda^G(X)$.
These polytopes are studied from symplectic point of view (see \cite{vergne:polmomentsurvey}) or
from an algebro-geometric point of view in \cite{Br:genface,man}.\\ 

In Section~\ref{sec:fprod}, we fix notation about parabolic fiber products and prove a useful result on their
linearized Picard group. Section~\ref{sec:HM} is concerned by the Hilbert-Mumford's numerical criterion of
 semistability.
In Section~\ref{sec:BB}, we recall some useful results about the Bialinicki-Birula's cells.
 Section~\ref{sec:Luna} recalls a useful consequence of Luna's Slice Etale Theorem.
In Section~\ref{sec:Cov}, we introduce the notion of well covering pair and prove their first relations with 
the $G$-cones.
In Section~\ref{sec:gen}, our general results about the faces of the $G$-cones are obtained.
In Section~\ref{sec:BeSj}, we apply our results to the cone $\lr(G,\hG)$ essentially by making more explicit 
the notion of well covering pair.\\

\noindent{\bf Convention.}
The ground field $\kk$ is assumed to be algebraic closed of characteristic zero.
The notation introduced in the environments ``{\bf Notation.}'' are fixed for all the sequence 
of the article.

\section{Preliminaries on parabolic fiber products}
\label{sec:fprod}

\begin{nota}
Let $\kk^*$ denote the multiplicative group of $\kk$.
  If $G$ is an affine algebraic group, $\Chi(G)$ denotes the group of  characters of $G$;
that is, of algebraic group homomorphisms from $G$ on $\kk^*$.
If $G$ acts algebraically  on a variety $X$, $X$ is said to be a {\it $G$-variety}.
As in \cite{GIT}, we denote by  $\Pic^G(X)$ the group of $G$-linearized line bundles on $X$.
If $\Li\in\Pic^G(X)$, ${\rm H}^0(X,\Li)$ denotes the $G$-module of regular sections of $\Li$.
\end{nota}

In this section we collect some useful properties of the fiber product.
Let us fix a reductive group $G$ and a parabolic subgroup $P$ of $G$.

\subsection{Construction} 
\label{sec:PicGP}

Let $Y$ be a $P$-variety.
Consider over $G\times Y$ the action of $G\times P$ given by the formula (with obvious notation):
$$
(g,p).(g',y)=(gg'p^{-1},py).
$$
Since the quotient map $G\longto G/P$ is a Zariski-locally trivial principal $P$-bundle; 
one can easily construct a quotient $G\times_PY$ of $G\times Y$ by the action of $\{e\}\times P$.
The action of $G\times\{e\}$ induces an action of $G$ on $G\times_PY$.
Moreover, the first projection $G\times Y\longto G$ induces a $G$-equivariant map  
$G\times_P Y\longto G/P$ which is a locally trivial fibration with fiber $Y$.

The class of a pair $(g,y)\in G\times Y$ in $\GPY$ is denoted by $[g:y]$.
If $Y$ is a $P$-stable locally closed subvariety of a $G$-variety $X$, it is well known that
the map
$$
\begin{array}{ccc}
\GPY&\longto& G/P\times X \\
\left[ g:y\right] & \longmapsto& (gP,gy)
\end{array}
$$ 
is an immersion; its image is the set of the $(gP,x)\in G/P\times X$ such that $g^{-1}x\in Y$.\\

Let $\nu$ be a character of $P$. If $Y$ is the field $\kk$ endowed with the action of $P$ defined by
$p.\tau=\nu(p^{-1})\tau$ for all $\tau\in\kk$ and $p\in P$, $G\times_P Y$ is a $G$-linearized 
line bundle on $G/P$.
We denote  by $\Li_\nu$ this element of $\Pic^G(G/P)$. 
In fact, the map $\Chi(P)\longto\Pic^G(G/P),\,\nu\longmapsto\Li_\nu$ is an isomorphism.

Let $B$ be a Borel subgroup of $G$ contained in $P$, and $T$ be a maximal torus contained in $B$.
Then,  $\Chi(P)$ identifies with a subgroup of $\Chi(T)$ which contains dominant weights.
For $\nu\in\Chi(P)$, $\Li_\nu$ is generated by its sections if and only if it has non zero sections
if and only if $\nu$ is dominant.
Moreover, ${\rm H}^0(G/P,\Li_\nu)$ is the dual of the simple $G$-module of highest weight $\nu$.
For $\nu$ dominant, $\Li_\nu$ is ample if and only if $\nu$ cannot be extended to a subgroup of $G$
bigger than $P$.

\subsection{Line bundles}

We are now interested in  the $G$-linearized line bundles on $\GPY$.

\begin{lemma}
\label{lem:picGPY}   
With above notation, we have:
\begin{enumerate}
        \item \label{ass1:picGPY}
        The map $\Li\longmapsto G\times_P\Li$ defines a morphism
        $$
        e\ :\ \Pic^P(Y)\longto \Pic^G(G\times_PY).
        $$
        \item The map $\iota\ :\ Y\longto\GPY,\,y\longmapsto [e:y]$ is a $P$-equivariant
        immersion. We denote by $\iota^*\ :\ \Pic^G(\GPY)\longto\Pic^P(Y)$ the associated
        restriction homomorphism.
        \item \label{ass3:picGPY}
        The morphisms $e$ and $\iota^*$ are the inverse one of each other; in particular, 
        they are isomorphisms.
        \item For any $\Li\in \Pic^G(\GPY)$, the restriction map 
        from ${\rm H}^0(\GPY,\Li)$ to ${\rm H}^0(Y,\iota^*(\Li))$ induces 
        a linear isomorphism 
        $${\rm H}^0(\GPY,\Li)^G\simeq {\rm H}^0(Y,\iota^*(\Li))^P.$$
\end{enumerate}
\end{lemma}

\begin{proof}
Let $\Mi$ be a $P$-linearized line bundle on $Y$.
Since the natural map $G\times \Mi\longto G\times_P \Mi$ is a categorical quotient, we have the following 
commutative diagram: 
$$
\begin{diagram}[height=4ex]
G\times \Mi &\rTo&G\times_P \Mi\\
\dTo&&\dTo_p\\
G\times Y&\rTo&G\times_P Y.
\end{diagram}
$$
Since $G\longto G/P$ is locally trivial, the map $p$ endows $G\times_P \Mi$ with a structure of 
line bundle on $G\times_P Y$.
Moreover, the action of $G$ on $G\times_P \Mi$ endows this line bundle with a $G$-linearization.
This proves Assertion~\ref{ass1:picGPY}. The second one is obvious.

By construction, the restriction of  $G\times_P \Mi$ to $Y$ is $\Mi$. 
So, $\iota^*\circ e$ is the identity map.
Conversely, let $\Li\in\Pic^G(\GPY)$. Then, we have:
$$
e\circ\iota^*(\Li)\simeq\{(gP,l)\in G/P\times \Li\::\:g^{-1}l\in\Li_{|Y}\}.
$$
The second projection induces an isomorphism from $e\circ\iota^*(\Li)$ onto $\Li$.
This ends the proof of Assertion~\ref{ass3:picGPY}.

The map ${\rm H}^0(\GPY,\Li)^G\longto {\rm H}^0(Y,\iota^*(\Li))^P$ is clearly well defined and injective.
Let us prove the surjectivity.
Let $\sigma\in{\rm H}^0(Y,\iota^*(\Li))^P$.
Consider the morphism
$$
\begin{array}{cccc}
\hat{\sigma}\::&G\times Y&\longto&G\times_P\Li\\
&(g,y)&\longmapsto&[g:\sigma(y)].
\end{array}
$$ 
Since $\sigma$ is $P$-invariant, so is $\hat{\sigma}$; and
$\hat{\sigma}$ induces a section of $G\times_P\Li$ over $G\times_PY$ which is $G$-invariant
and extends $\sigma$.
\end{proof}

\section{Numerical criterion of Hilbert-Mumford}
\label{sec:HM}

\begin{nota}
  If $x$ is a point of a $G$-variety, we will denote by $G_x$ its isotropy subgroup and by
$G.x$ its orbit.
\end{nota}

We will use classical results in Geometric Invariant Theory (GIT)  
about the numerical criterion of Hilbert-Mumford. In this section, we present
these results and give some useful complements.
Let us fix a connected reductive group $G$ acting on an irreducible projective 
algebraic variety $X$.

\subsection{Convex geometry}
\label{sec:convex}

In this subsection, we recall some useful facts about polytopes.
Let $E$ be an Euclidean space whose the scalar product is denoted by 
$\langle\cdot\,,\,\cdot\rangle$ and norm $\Vert\cdot\Vert$. 
A {\it polytope} in $E$ is the convex hull of finitely many points.
Let $\P$ be a polytope which does not contain $0$.
The real $\inf_{x\in\P}\Vert x\Vert$ is called the {\it distance from $0$
to $\P$} and denoted by $d(O,\P)$.
The orthogonal projection theorem asserts that there exists a unique 
$x_0\in\P$ such that $d(0,\P)=\Vert x_0\Vert$. Moreover,
$$
d(0,\P)=\sup_{y\in E-\{0\}}\inf_{x\in\P}\frac{\langle y,x\rangle}{\Vert y\Vert},
$$
and $x_0'=\frac{x_0}{\Vert x_0\Vert}$ is the unique vector with norm 1 such that
$$
d(0,\P)=\inf_{x\in\P}\langle x_0',x\rangle.
$$
In particular, the set of $x\in\P$ such that $\langle x_0',x\rangle=d(0,\P)$ is a face of $\P$; 
we call it {\it the face viewed from $0$} and denote it by $\Face(0)$.

We will use the following easy fact:

\begin{lemma}
\label{lem:convex}
Let $r$ be a positive integer.
We consider the map which associate to any element $\underline{x}=(x_1,\cdots,x_r)$ of $E^r$  
their convex hull ${\rm Conv}(x_1,\cdots,x_r)=\P(\underline{x})$.

Then, the set $U$ of $\underline{x}\in E^r$ such that $0\not\in \P(\underline{x})$ is open in $E^r$. 
To any such $\underline{x}$, we associate the set $I(\underline{x})$ of the $i$'s such that
$x_i$ belongs to the face of $\P(\underline{x})$ viewed from $0$. 
Then, there exists an open neighborhood $U$ of $\underline{x}$ such that for any 
$\underline{x}'\in U$ we have $I(\underline{x}')\subset I(\underline{x})$.
\end{lemma}

\begin{proof}
  The only point is that the application which maps $\underline{x}$ to the orthogonal projection
of $0$ on  $\P(\underline{x})$ is continuous.
\end{proof}

\subsection{An Ad Hoc notion of semistability}

As in the introduction, for any  $G$-linearized line bundle $\Li$ on $X$, we consider the
following set of {\it semistable points}:
$$
X^{\rm ss}(\Li)=\left \{
x\in X \mbox{ : }
\exists n>0\mbox{ and }
\sigma\in{\rm H}^0(X,\Li^{\otimes n})^G \mbox{ such that }\sigma(x)\neq 0
\right \}.
$$ 
To precise the acting group, we sometimes denote $X^{\rm ss}(\Li)$ by $X^{\rm ss}(\Li,G)$.

The subset $X^{\rm ss}(\Li)$ is open and stable by $G$.
A point $x$ which is not semistable is said to be {\it unstable}; and, we set
$X^{\rm us}(\Li)=X-X^{\rm ss}(\Li)$.\\

\begin{remark}
  Note that this definition of $X^{\rm ss}(\Li)$ is NOT standard. Indeed, one usually imposes that the open
subset defined by the non vanishing of $\sigma$ to be affine. This property which is useful to 
construct a good quotient is automatic if $\Li$ is ample but not in general; hence, 
our definition coincides with the usual one if $\Li$ is ample.
\end{remark}\\

If $\Li$ is ample, there exists a categorical quotient:
$$
\pi\::\:X^{\rm ss}(\Li)\longto X^{\rm ss}(\Li)\quot G,
$$
such that $X^{\rm ss}(\Li)\quot G$ is a projective variety and $\pi$ is affine.
A point $x\in X^{\rm ss}(\Li)$ is said to be {\it stable} if $G_x$ is finite and $G.x$ is 
closed in $X^{\rm ss}(\Li)$. 
Then, for all stable point $x$ we have $\pi^{-1}(\pi(x))=G.x$; and the set 
$X^{\rm s}(\Li)$ of stable points is open in $X$.

The following lemma is easy and well known. It will be very useful here.

\begin{lemma}
\label{lem:LGx}
Let $\Li$ be an $G$-linearized line bundle on $X$
and  $x\in X$ be a point {\it semistable} for $\Li$.

Then, the restriction of $\Li$ to $G.x$ is of finite order.
\end{lemma}

\begin{proof}
Let us recall that for any $\Li\in\Pic^G(G.x)$, the action of $G_x$ on the fiber over $x$ in $\Li$ 
determines a character $\mu^\Li(x,G_x)$ of $G_x$. Moreover, the map $\Li\mapsto\mu^\Li(x,G_x)$ is an injective
homomorphism.

Now, let $\Li$ be a $G$-linearized line bundle on $X$ such that the character $\mu^\Li(x,G_x)$ is of infinite order.
It remains to prove that $x$ is unstable for $\Li$.
Let $\sigma$ be a $G$-invariant section of $\Li^{\otimes n}$ for some $n>0$.
Then $\sigma(x)$ is a $G_x$ fix point of the fiber in $\Li^{\otimes n}$ over $x$.
Since, $n.\mu^\Li(x,G_x)$ is non trivial, $\sigma(x)$ must be zero.
So, $x$ is unstable.
\end{proof}\\

\subsection{The functions $\mu^\bullet(x,\lambda)$}
\label{sec:mu}

Let $\Li\in \Pic^G(X)$.
Let $x$ be a point in $X$ and $\lambda$ be a one parameter subgroup of
$G$.
Since $X$ is complete,  $\lim_{t \to 0}\lambda(t)x$ exists; let $z$
denote this limit.
The image of $\lambda$ fixes $z$ and so the group $\kk^*$ acts
via $\lambda$ on the fiber $\Li_{z}$.
This action defines a character of $\kk^*$, that is, an element of
$\ZZ$ denoted by $\mu^\Li(x,\lambda)$. 
One can immediately prove that the numbers $\mu^\Li(x,\lambda)$ 
satisfy the following properties:
\begin{enumerate}
\item $\mu^\Li(g\cdot x,g\cdot \lambda\cdot g^{-1})=\mu^\Li(x,\lambda)$ 
for any $g\in G$;
\item  the map $\Li\mapsto \mu^\Li(x,\lambda)$
  is a homomorphism from $\Pic^G(X)$ to $\ZZ$;
\item \label{prop*mu} for any $G$-variety $Y$ and for any $G$-equivariant morphism $f\::\:X\longto Y$, 
$\mu^{f^*(\Li)}(x,\lambda)=\mu^\Li(f(x),\lambda)$, where $x\in X$ and $\Li$ is a $G$-linearized line bundle on $Y$. 
\end{enumerate}

\noindent
A less direct property of the function $\mu^\Li(x,\lambda)$ is 
\begin{prop}
\label{prop:sgnmu}
Let $\Li,\,x,\,\lambda$ and $z$ be as above.
Let $\tilde{x}$ be a non zero point in the fiber in $\Li$ over $x$.
We embed $X$ in $\Li$ by the zero section.
Then, we have

\begin{enumerate}
\item if $\mu^\Li(x,\lambda)>0$, then $\lambda(t)\tilde{x}$ tends to $z$ when $t\to 0$;
\item if $\mu^\Li(x,\lambda)=0$, then $\lambda(t)\tilde{x}$ tends to a non zero point $\tilde{z}$
in the fiber in $\Li$ over $z$ when $t\to 0$;
\item if $\mu^\Li(x,\lambda)<0$, then $\lambda(t)\tilde{x}$ has no limit in $\Li$ when $t\to 0$.
\end{enumerate}
\end{prop}

\begin{proof}
Set $V=\{\lambda(t).x\:|\;t\in\kk^*\}\cup \{z\}$: it is a locally closed subvariety of $X$ stable by the 
action of $\kk^*$ via $\lambda$.
Moreover, $z$ is the  unique closed orbit of $\kk^*$ in $V$.
So, \cite[Lemma~7]{GeomDedic}  implies that $\Pic^{\kk^*}(V)$ is isomorphic to $\Chi(\kk^*)$; 
and finally that for all $\Li\in\Pic^G(X)$ the restriction $\Li$ to $V$ is the trivial line bundle
endowed with the action of $\kk^*$ given by $\mu^\Li(x,\lambda)$. 
The proposition follows immediately.
\end{proof}\\

The numbers $\mu^\Li(x,\lambda)$ are used in \cite{GIT} to give a 
numerical criterion for stability with respect to an ample $G-$linearized
line bundle $\Li$:
$$
\begin{array}{l}
  x\in X^{\rm ss}(\Li)\iff \mu^\Li(x,\lambda)\leq 0 \mbox{ for all
    one parameter subgroup $\lambda$},\\
x\in X^{\rm s}(\Li)\iff \mu^\Li(x,\lambda)< 0 \mbox{ for all non trivial 
$\lambda$}.
\end{array}
$$

A line bundle $\Li$ over $X$ is said to be {\it semiample} if a positive power of $\Li$ is base point free.
With our notion of semistability, Hilbert-Mumford's theorem admits the following direct generalization:

\begin{lemma}
\label{lem:HMsemiample}
Let $\Li$ be a $G$-linearized  line bundle over $X$ and $x$ a point in $X$.
Then,
\begin{enumerate}
\item if $x$ is semistable for $\Li$, $\mu^\Li(x,\lambda)\leq 0$ for any one parameter subgroup $\lambda$ of $G$;
\item\label{assGeomDedic} for a one parameter subgroup $\lambda$ of $G$, 
if $x$ is semistable for $\Li$ and  $\mu^\Li(x,\lambda)=0$, then $\lim_{t\to 0}\lambda(t)x$ is semistable for $\Li$;
\item \label{ass:HMsa} if in addition $\Li$ is semiample,  $x$ is semistable for $\Li$ if and only if  $\mu^\Li(x,\lambda)\leq 0$ 
for any one parameter subgroup $\lambda$ of $G$.
\end{enumerate}
\end{lemma}

\begin{proof}
Assume that $x$ is semistable for $\Li$ and consider a $G$-invariant section $\sigma$ of $\Li^{\otimes n}$ which does
not vanish at $x$.
Since $\lambda(t)\sigma(x)=\sigma(\lambda(t)x)$ tends to $\sigma(z)$  when $t\to 0$, Proposition~\ref{prop:sgnmu} shows that 
$\mu^\Li(x,\lambda)\leq 0$. 
If in addition $\mu^\Li(x,\lambda)=0$, Proposition~\ref{prop:sgnmu} shows that $\sigma(z)$ is non zero; and so that 
$z$ is semistable for $\Li$. This proves the two first assertions. 

Assume now that $\Li$ is semiample.
Let $n$ be a positive integer such that $\Li^{\otimes n}$ is base point free.
Let $V$ denote the dual of the  space of global sections of $\Li^{\otimes n}$: 
$V$ is a finite dimensional $G$-module. Moreover, the usual map
$\phi\,:\,X\longto \PP(V)$ is $G$-equivariant. Let $Y$ denote the image of $\phi$
and $\Mi$ denote the restriction of $\Oun$ to $Y$.

Then $\Li^{\otimes n}$ is the pullback of $\Mi$ by $\phi$.
Since $X$ is projective, $\phi$ induces isomorphisms from 
${\rm  H}^0(Y,\Mi^{\otimes k})$ onto  ${\rm H}^0(X,\Mi^{\otimes nk})$ (for all $k$).
So, $X^{\rm ss}(\Li)=\phi^{-1}(Y^{\rm ss}(\Mi))$.
We deduce the last assertion of the lemma by applying  Hilbert-Mumford's criterion to $Y$ and $\Mi$ and 
Property~\ref{prop*mu} of the functions $\mu^\bullet(x,\lambda)$.
\end{proof}

\bigskip
\begin{remark}
  \begin{enumerate}
  \item
  If $\Li$ is ample, Assertion~\ref{assGeomDedic} of Lemma~\ref{lem:HMsemiample} is Lemma~3 in \cite{GeomDedic}. 
\item The proof of Assertion~\ref{ass:HMsa} shows that a lot of properties of semistability for an ample line bundle
are also available  for semiample line bundles (see Propositions~\ref{PropMtore} and \ref{prop:Mfini}, and
Theorems~\ref{thKempf} and \ref{thNess} below).
 \end{enumerate}
\end{remark}

\subsection{Definition of the functions $\mbox{M}^\bullet(x)$}

\begin{nota}
Let $\Gamma$ be any affine algebraic group.
The neutral component of $\Gamma$ is denoted by $\Gamma^\circ$.
 Let $\Y(\Gamma)$ denote the set of one parameter subgroups of $\Gamma$; that is, 
of group  homomorphisms from $\kk^*$ in $\Gamma$.
Note that if $\Gamma^\circ$ is a torus, $\Y(\Gamma)$ is a group.

If $\Lambda$ is an abelian group, we denote by $\Lambda_\QQ$ (resp. $\Lambda_\RR$) the tensor product of 
$\Lambda$ with $\QQ$ (resp. $\RR$) over $\ZZ$.
\end{nota}

Let $T$ be a maximal torus of $G$.
The Weyl group $W$ of $T$ acts linearly on $\Y(T)_\RR$.
Since $W$ is finite,
there exists a $W$-invariant Euclidean norm (defined over $\QQ$) 
$\Vert\cdot \Vert$ on $\Y(T)_\RR$. 
On the other hand, if $\lambda\in \Y(G)$ there exists $g\in G$ such
that  $g\cdot \lambda\cdot g^{-1}\in\Y(T)$. 
Moreover, if two elements of $\Y(T)$
are conjugate by an element of $G$, then they are by an element of the
normalizer of $T$ (see  \cite[Lemma~2.8]{GIT}).
This allows  to define the norm of $\lambda$ by
$\Vert\lambda\Vert=\Vert g\cdot \lambda \cdot g^{-1}\Vert$. 

Let $\Li\in \Pic^G(X)$. One can now introduce the following notation:
$$
\begin{array}{cc}
\overline{\mu}^\Li(x,\lambda)=\frac{\mu^\Li(x,\lambda)}{\Vert
  \lambda\Vert},\hspace{0.2cm}&
\displaystyle{ \mbox{M}^\Li(x)=
\sup_{\lambda\in\Y(G)}\overline{\mu}^\Li(x,\lambda).}
\end{array}
$$

In fact, we will see in Corollary~\ref{cor:Mfini} that  $\mbox{M}^\Li(x)$ is finite.

\subsection{$\mbox{M}^\bullet(x)$ for a torus action}
\label{Mtore}

\begin{nota}
If $Y$ is a variety, and $Z$ is a part of $Y$, the closure of $Z$ in $Y$ will be denote by $\overline{Z}$.
If $\Gamma$ is a group acting on $Y$, $Y^\Gamma$ denote the set of fix point of $\Gamma$ in $Y$.

If $V$ is a finite dimensional vector space, and $v$ is a non zero vector in $V$, $[v]$ denote the class of
$v$ in the projective space $\PP(V)$.
If $V$ is a $\Gamma$-module, and $\chi$ is a character of $\Gamma$, we denote by $V_\chi$ the set of $v\in V$
such that $g.v=\chi(g)v$ for all $g\in\Gamma$. 
\end{nota}

In this subsection we assume that $G=T$ is a torus.
Let $z$ be a point of $X$ fixed by $T$.
The action of $T$ on the fiber $\Li_z$ over the point $z$ in the $T$-linearized 
line bundle $\Li$ define a character $\mu^\Li(z,T)$ of $T$; we obtain a morphism
$$
\mu^\bullet(z,T)\::\:\Pic^T(X)\longto \Chi(T).
$$
For any point $x$ in $X$, we denote by $\P_T^\Li(x)$ the convex hull in $\Chi(T)_\RR$ of the characters
$-\mu^\Li(z,T)$ for $z\in\overline{T.x}^T$.

The following proposition is an adaptation of a result of L. Ness and 
gives a pleasant interpretation of the number $\mbox{M}^\Li(x)$:

\begin{prop}
\label{PropMtore}
Let $\Li$ be a semiample $T$-linearized line bundle on $X$.
With the above notation, we have:
\begin{enumerate}
\item
The point $x$ is unstable if and only if $0$ does not belong to $\P_T^\Li(x)$.
In this case, ${\rm M}^\Li(x)$ is the distance from $0$ to $\P_T^\Li(x)$.
\item
If $x$ is semistable, the opposite of ${\rm M}^\Li(x)$ is the 
distance from $0$ to the  boundary of $\P_T^\Li(x)$.
\item
There exists $\lambda\in\Y(T)$ such that
$\overline{\mu}^\Li(x,\lambda)={\rm M}^\Li(x)$.
If moreover $\lambda$ is indivisible, we call it an 
{\it adapted one parameter subgroup for $x$}.
\item
If $x$ is unstable, there exists a unique adapted 
one parameter subgroup for $x$.
\end{enumerate}
\end{prop}

\begin{proof}
Since $\Li$ is semiample, there exist a positive integer $n$,  a $T$-module $V$,
and a $T$-equivariant morphism $\phi\::\:X\longto \PP(V)$
such that $\Li^{\otimes n}=\phi^*({\cal O}(1))$.
Since $\mu^\bullet(z,T)$ is a morphism, we have: $\P_T^{\Li^{\otimes n}}(x)=n\P_T^\Li(x)$ for all $x$.
Moreover, $\overline{\mu}^{\Li^{\otimes n}}(x,\lambda)=n\overline{\mu}^{\Li}(x,\lambda)$, for all $x$ and $\lambda$;
so, $M^{\Li^{\otimes n}}(x)=nM^\Li(x)$.
As a consequence, it is sufficient to prove the proposition for $\Li^{\otimes n}$; in other words, we may assume that $n=1$.

Let us recall that:
$$
V=\bigoplus_{\mbox{\scriptsize $\chi\in \Chi(T)$}} V_{\chi}.
$$

Let $x\in X$ and $v\in V$ such that $[v]=\phi(x)$.
There exist unique vectors $v_\chi\in V_\chi$ such that 
$v=\sum_{\chi}v_\chi$.
Let $\Q$ be the convex hull in $\Chi(T)_\RR$ of the $\chi$'s such that $v_\chi\neq 0$.
It is well known (see~\cite{Oda}) that the fixed point of $T$ in $\overline{T.[v]}$ are exactly the
$[v_\chi]$'s with $\chi$ vertex of $\Q$. 
Moreover, $T$ acts by the character $-\chi$ on the fiber ${\cal O}(1)_{[v_\chi]}$ over 
$[v_\chi]$ in ${\cal O}(1)$.
One  deduces that $\Q= \P_T^{\Li}(x)$. 
So, it is sufficient to prove the proposition with $Q$ in place of 
$ \P_T^{\Li}(x)$ and $\phi(x)=x$ (because of our non-standard definition of
semistability): this is a statement of \cite{Ne1}.
\end{proof}

\subsection{Properties of $M^\bullet(x)$}
\label{Mfini}

\begin{nota}
We will denote by $\Pic^G(X)^+$ (resp.   $\Pic^G(X)^{++}$) the set of semiample (resp. ample) $G$-linearized
line bundles on $X$. 
\end{nota}

An indivisible one parameter subgroup $\lambda$ of $G$ is said to be
{\it   adapted for   $x$ and $\Li$} if and only if
$\overline{\mu}^\Li(x,\lambda)=\mbox{M}^\Li(x)$.  
Denote by $\Lambda^\Li(x)$ the set of adapted one parameter subgroups
for $x$.

Using the fact that any one parameter subgroup is conjugated to one in a given
torus, Proposition~\ref{PropMtore} implies easily 
 
\begin{coro}
\label{cor:Mfini}
\begin{enumerate}
\item The numbers $\mbox{M}^\Li(x)$ are finite (even if $\Li$ is not semiample,
see  Proposition 1.1.6 in \cite{DH}).
\item If $\Li$ is semiample, $\Lambda^\Li(x)$ is not empty.
\end{enumerate}
\end{coro}

Now, we can reformulate the  numerical criterion for stability: if 
$\Li$ is semiample, we have
$$
\begin{array}{cc}
X^{\rm ss}(\Li)=\{x\in X : \mbox{M}^\Li(x)\leq 0\},&
X^{\rm s}(\Li)=\{x\in X : \mbox{M}^\Li(x)< 0\}.
\end{array}
$$

The following proposition is a result of finiteness for the set of functions $M^\bullet(x)$.
It will be used to understand how $X^{\rm ss}(\Li)$ depends on $\Li$ (see Lemma~\ref{lem:Xcircmono} below).

 \begin{prop}
\label{prop:Mfini}
When $x$ varies in $X$, one obtains only a finite number of functions
$M^\bullet(x)\ :\ \Pic^G(X)^+\longto\RR$.
\end{prop}

\begin{proof}
Let $T$ be a maximal torus of $G$.
Consider the partial forgetful map $r^T\::\:\Pic^G(X)\longto\Pic^T(X)$.
Since $M^\bullet(x)=\max_{g\in G}M^{r^T(\bullet)}(g.x)$, it is sufficient to prove
the proposition for the torus $T$.

If $z$ and $z'$ belong to the same irreducible component $C$ of $X^T$, 
the morphisms $\mu^\bullet(z,T)$ and $\mu^\bullet(z',T)$ are equal:
we denote by $\mu^\bullet(C,T)$ this morphism.
 
By Proposition~\ref{PropMtore}, $M^\Li(x)$ only depends on $\P_T^\Li(x)$, which only depends
on the set of irreducible components of $X^T$ which intersects $\overline{T.x}$.
Since, $X^T$ has finitely many irreducible components, the proposition follows.
\end{proof}\\

\begin{remark}
Proposition~\ref{prop:Mfini} implies that the open subsets of $X$ which can be realized as 
$X^{\rm ss}(\Li)$ for some semiample $G$-linearized line bundle $\Li$ on $X$ is finite.
This is  a result of Dolgachev and Hu (see Theorem~3.9 in~\cite{DH}; see also \cite{Schmitt}).
\end{remark}

\subsection{Adapted one parameter subgroups}
\label{sec:aops}

To describe $\Lambda^\Li(x)$, we need some additional notation.
To the one parameter subgroup $\lambda$ of $G$, we associate the
parabolic subgroup (see \cite{GIT}):
$$
P(\lambda)=\left \{
g\in G \::\:
\lim_{t\to 0}\lambda(t).g.\lambda(t)^{-1} 
\mbox{  exists in } G \right \}.
$$
The unipotent radical of $P(\lambda)$ is 
$$
U(\lambda)=\left \{
g\in G \::\:
\lim_{t\to 0}\lambda(t).g.\lambda(t)^{-1}=e 
\right \}.
$$
Moreover, the centralizer $G^\lambda$ of the image of $\lambda$ in $G$ is a Levi subgroup
of $P(\lambda)$.
For $p\in P(\lambda)$, we set $\overline{p}=\lim_{t\to 0}\lambda(t).p.\lambda(t)^{-1}$.
Then, we have the following short exact sequence:
\begin{diagram}
1&\rTo&U(\lambda)&\rTo&P(\lambda)&\rTo^{p\mapsto\overline{p}}&G^\lambda&\rTo&1.
\end{diagram} 

For $g\in P(\lambda)$, we have
$\mu^\Li(x,\lambda)=\mu^\Li(x,g\cdot\lambda\cdot g^{-1})$. 
The following theorem due to G. Kempf is a generalization of  the last
assertion of Proposition~\ref{PropMtore}.

\begin{theo}(see \cite{Ke})
\label{thKempf}
Let $x$ be an unstable point for a semiample $G$-linearized line
bundle $\Li$. Then:
\begin{enumerate}
\item All the $P(\lambda)$ for $\lambda\in\Lambda^\Li(x)$ are equal. We
  denote by $P^\Li(x)$ this subgroup.
\item Any two elements of $\Lambda^\Li(x)$ are conjugate by an element
  of $P^\Li(x)$.   
\end{enumerate}
\end{theo}

We will also use the following theorem of L. Ness.

\begin{theo}(Theorem 9.3 in \cite{Ne2})
\label{thNess}
Let $x$ and $\Li$ be as in the above theorem.
Let $\lambda$ be an adapted one parameter subgroup for $x$ and $\Li$.
We consider $y=\lim_{t\to 0}\lambda(t)\cdot x$.
Then, $\lambda\in\Lambda^\Li(y)$ and ${\rm M}^\Li(x)={\rm M}^\Li(y)$.
\end{theo}

\subsection{Stratification of $X$ induced from $\Li$}
\label{sec:strat}

Let $\Li$ be an ample $G$-linearized line bundle on $X$.
If $d>0$ and $\langle\tau\rangle$ is a conjugacy class of one parameter subgroups of
$G$, we set:
$$
S^\Li_{d,\langle\tau\rangle}=\left \{ 
x\in X \::\: \mbox{M}^\Li(x)=d \mbox{ and  } 
\Lambda^\Li(x)\cap\langle\tau \rangle \neq\emptyset
\right \}.
$$
If ${\cal T}$ is the set of conjugacy classes of one parameter subgroups,
the previous section gives us the following partition of $X$:
\begin{eqnarray}
  \label{eq:Hess}
X=X^{\rm ss}(\Li)\cup\bigcup_{d>0,\;\langle\tau \rangle \in{\cal T}}S^\Li_{d,\langle\tau \rangle }.
\end{eqnarray}
W. Hesselink showed in \cite{He} that this union is a finite
stratification by $G$-stable locally closed subvarieties of $X$.
We will call it the {\it stratification induced from $\Li$}. 

To describe the geometry of these stratum, we need additional notation.
For $\lambda\in \langle\tau \rangle $, we set:
$$
S^\Li_{d,\lambda}:=\{x\in S^\Li_{d,\langle\tau \rangle }\ :\ \lambda\in\Lambda^\Li(x)\},
$$
and 
$$
Z^\Li_{d,\lambda}:=\{x\in S^\Li_{d,\lambda}\ :\ \lambda(\kk^*){\rm\ fixes\ }x\}.
$$
By Theorem~\ref{thNess}, we have the map
$$
\Lim\ :\ S^\Li_{d,\lambda}\longto Z^\Li_{d,\lambda},\ x\longmapsto\lim_{t\to 0} \lambda(t).x.
$$
The proof of the following result can be found in \cite[Section~1.3]{Ki}.

\begin{prop}
\label{prop:Ki}
With above notation, if $d$ is positive, we have:
\begin{enumerate}
  \item $Z^\Li_{d,\lambda}$ is open in $X^\lambda$ and stable by $G^\lambda$;
  \item $S^\Li_{d,\lambda}=\{x\in X\ :\  \lim_{t\to 0} \lambda(t).x\in Z^\Li_{d,\lambda}\}$ and is
stable by $P(\lambda)$;
 \item \label{ass4:Ki}
 there is a birational morphism $G\times_{P(\lambda)} S^\Li_{d,\lambda}\longto S^\Li_{d,\langle\lambda \rangle }$,
 which is an isomorphism if $S^\Li_{d,\langle\lambda \rangle }$ is normal.
\end{enumerate}
\end{prop}

\subsection{A description of $Z^\Li_{d,\lambda}$}
\label{sec:Hlambda}

Let $\lambda$ be a one parameter subgroup of $G$.
Let $Z$ denote the neutral component of the center of $G^\lambda$ and 
$G^{\rm ss}$ be the maximal semisimple subgroup  of $G^\lambda$.
The product induces an isogeny $Z\times G^{\rm ss}\longto G^\lambda$.
Let $T_1$ be a maximal torus of $G^{\rm ss}$. Set $T=Z.T_1$. Note that $T$  is a 
maximal torus of $G^\lambda$ and $G$.
 Let $S$ be the subtorus of $Z$ such that $Y(S.T_1)_\RR$ is the hyperplane of $Y(T)_\RR$
orthogonal to $\lambda$. Set $H^\lambda=S.G^{\rm ss}$.
The map $\GG_m\times H^\lambda\longto G^\lambda,(t,h)\longmapsto\lambda(t)h$ is an isogeny.

\begin{theo}[Ness-Kirwan]
\label{th:NessKi}
Let $\Li$ be an ample $G$-linearized line bundle on $X$. The one parameter subgroup $\lambda$ is assumed to
be indivisible.
  Let $x\in X^\lambda$ be such that $\mu^\Li(x,\lambda)>0$.

Then, $\lambda$ is adapted to $x$ and $\Li$ if and only if 
$x$ is semistable for $\Li$ and the action of $H^\lambda$.
\end{theo}

Theorem~\ref{th:NessKi} is a version of  \cite[Theorem~9.4]{Ne2}. 
Whereas the Ness' proof works without changing, the statement in \cite{Ne2} is not correct.
In \cite[Remark~12.21]{Ki}, F. Kirwan made the above correction.

\subsection{The open stratum}

Let $\Li$ be an ample $G$-linearized line bundle on $X$. 
We denote by $X^\circ(\Li)$ the open stratum of Stratification~\ref{eq:Hess}.
We have the following characterization of $X^\circ(\Li)$:

\begin{lemma}
\label{lem:dmin}
If $X^{\rm ss}(\Li)\neq\emptyset$, set $d_0=0$; else, set $d_0=\min_{x\in X}M^\Li(x)$. 
Then,   $X^\circ(\Li)$ is the set of $x\in X$ such that $M^\Li(x)\leq d_0$.
\end{lemma}

\begin{proof}
If $d_0=0$,   $M^\Li(x)\leq 0$ if and only if $x\in X^{\rm ss}(\Li)=X^\circ(\Li)$.
We now assume that $d_0>0$.

Up to changing $\Li$ by a positive power, one may assume that 
there exists a $G$-module $V$ such that $X$ is contained in $\PP(V)$ and 
$\Li$ is the restriction of the $G$-linearized line bundle $\O(1)$.

Let us fix the positive real $d$ and a one parameter subgroup $\lambda$ such that 
$X^\circ(\Li)=S^\Li_{d,\langle\lambda\rangle}$.
For $i\in \ZZ$, set $V_i=\{v \in V\,|\,\lambda(t)v=t^iv\}$.
Set $V^+=\oplus_{i> d \Vert\lambda\Vert}$, $C=\PP(V_{d \Vert\lambda\Vert})\cap X$
and $\overline{C^+}=\PP(V_{d\Vert\lambda\Vert}\oplus V^+)\cap X$.

We first prove that $d=d_0$; that is, that $M^\Li(x)\geq d$ for all $x\in X$.  
Consider  the morphism 
$\eta\,:\,G\times_{P(\lambda)}\overline{C^+}\longto X$.
Since $G/P(\lambda)$ is projective, $\eta$ is proper; but, the image of $\eta$ contains
$S^\Li_{d,\langle\lambda\rangle}$; so, $\eta$ is surjective.
Let $x\in X$. There exists $g\in G$ such that $gx\in  \overline{C^+}$.
Then, $M^\Li(x)\geq \bar\mu^\Li(gx,\lambda)\geq d$.

Conversely, let $x\in X$ be such that $M^\Li(x)=d_0$.
 Let $g\in G$ such that $gx\in\overline{C}^+$.
 Let $v_1\in V_{d\Vert\lambda\Vert}$ and $v_2\in V^+$ such that $gx=[v_1+v_2]$.
Since $\overline{\mu}(gx,\lambda)\leq d_0$, $v_1$ is non zero; 
so, $\overline\mu(gx,\lambda)=d_0=M^\Li(x)$. In particular, $\lambda$ is adapted to $x$ and $\Li$; 
that is, $x\in X^\circ(\Li)$.
\end{proof}\\

We will need the following result of monotonicity for the function $\Li\mapsto X^\circ(\Li)$:

\begin{lemma}
  \label{lem:Xcircmono}
With above notation, there exists an open neighborhood $U$ of $\Li$ in $\Lambda_\QQ^{++}$ 
such that for all $\Li'\in U$, $X^\circ(\Li')\subset X^\circ(\Li)$.
\end{lemma}

\begin{proof}
By Proposition~\ref{prop:Mfini} and Lemma~\ref{lem:dmin}, there exists only finitely many open subsets 
 of $X$ which are of the form  $X^\circ(\Mi)$ for some ample $\Mi$.
Let  $X^\circ_1,\cdots,X^\circ_s$ those which are not contained in $X^\circ(\Li)$.
For each $i$, fix $x_i\in  X^\circ_i-X^\circ(\Li)$.
It remains to prove that for each $i$, there exists $U_i$ such that $x_i\not\in X^\circ(\Li')$ for all 
$\Li'\in U_i$. Indeed, $U=\cap_iU_i$ will work.

Set $d_0$ be as in Lemma~\ref{lem:dmin}. Since $x_i\in X^\circ(\Li')$, $M^\Li(x)>d_0$.
Let $e$ be such that  $M^\Li(x)>e>d_0$.
By Proposition~\ref{prop:Mfini}, there exists an open  neighborhood $U'$ of $\Li$ such that 
$\min_{x\in X}M^{\Li'}(x)>e$ for all $\Li'\in U'$.
Moreover, there exists $U''$, such that $M^{\Li'}(x_i)<e$ for all $\Li'\in U'$.
Lemma~\ref{lem:dmin} implies that for all $\Li'\in U'\cap U''$, $x_i\not\in X^\circ(\Li')$.
\end{proof}

\section{Bialynicki-Birula cells}
\label{sec:BB}

\subsection{Bialynicki-Birula's theorem}

\begin{nota}
Let $X$ be a $G$-variety, $H$ be a subgroup of $G$ and $\lambda$ be a one parameter subgroup of $G$.
The centralizer of $H$ in $G$ (that is, the set of fix point of $H$ acting on $G$ by conjugacy) will 
be denoted by $G^H$. 
The set of fix points of the image of $\lambda$ will be denoted by $X^\lambda$; the centralizer of this image
will be denoted by $G^\lambda$.

If $Y$ is a locally closed subvariety of $X$, and $\Li$ is a line bundle on $X$, $\Li_{|Y}$ will denote 
the restriction of $\Li$ to $Y$.
\end{nota}

Let $X$ be a complete $G$-variety. Let $\lambda$ be a one parameter subgroup of $G$.
Let $C$ be an irreducible component of $X^\lambda$. 
Since $G^\lambda$ is connected, $C$ is a $G^\lambda$-stable closed subvariety of $X$.
We set:
$$
C^+:=\{x\in X\::\: \lim_{t\to 0}\lambda(t)x\in C\}.
$$
Then, $C^+$ is a locally closed subvariety of $X$ stable by $P(\lambda)$.
Moreover, the map $\Lim\::\:C^+\longto C,\,x\longmapsto\lim_{t\to 0}\lambda(t)x$ is 
a morphism satisfying:
$$
\forall (l,u)\in G^\lambda\times U(\lambda)\hspace{1cm}
\Lim(lu.x)=l\Lim(x).
$$

Let $x\in X^\lambda$. We consider the natural action of $\kk^*$ induced by $\lambda$
on the Zariski tangent space $T_xX$ of $X$ at $x$. We consider the following $\kk^*$-submodules
of $T_xX$:
$$
\begin{array}{l}
T_xX_{>0}=\{\xi\in T_xX\::\:\lim_{t\to 0}\lambda(t)\xi=0\},\\
T_xX_{<0}=\{\xi\in T_xX\::\:\lim_{t\to 0}\lambda(t^{-1})\xi=0\},\\
T_xX_0=(T_xX)^\lambda,\ T_xX_{\geq 0}=T_xX_{>0}\oplus T_xX_{0}
{\rm\ and\ }T_xX_{\leq 0}=T_xX_{<0}\oplus T_xX_{0}.
\end{array}
$$ 

A classical result of Bialynicki-Birula (see~\cite{BB}) is

\begin{theo}
\label{th:BB}
Assuming in addition that $X$ is smooth, we have:
\begin{enumerate}
\item $C$ is smooth and for all $x\in C$ we have $T_xC=T_xX_0$;
\item $C^+$ is smooth and irreducible and for all $x\in C$ we have $T_xC^+=T_xX_{\geq 0}$;
\item the morphism $\Lim\::\:C^+\longto C$ induces a structure of vector bundle on $C$ with fibers
isomorphic to $T_xX_{>0}$.
\end{enumerate}
\end{theo}

\subsection{Line bundles on $C^+$}

We will need some results about the line bundles on $C^+$.
Let $\Li$ be a $P(\lambda)$-linearized line bundle on $C^+$.
Since $\lim_{t\to 0}\lambda(t).x$ belongs to $C$, the number $\mu^\Li(x,\lambda)$ is well defined, for
$x\in C^+$.
Moreover, since $C$ is irreducible, $\mu^\Li(x,\lambda)$ does not depend on $x\in C^+$;
we denote by $\mu^\Li(C,\lambda)$ this integer.

\begin{prop}
\label{prop:picC+}
We assume that $X$ is smooth. Then, we have:
\begin{enumerate}
\item \label{ass1:picC+}
The restriction map $\Pic^{P(\lambda)}(C^+)\longto \Pic^{G^\lambda}(C)$ is an isomorphism. 

Let $\Li\in\Pic^{P(\lambda)}(C^+)$. 
\item \label{ass2:picC+}
If $\mu^\Li(C,\lambda)\neq 0$, ${\rm H}^0(C,\Li_{|C})^\lambda=\{0\}$.
\item \label{ass3:picC+}
If $\mu^\Li(C,\lambda)=0$, the restriction map induces an isomorphism
from ${\rm H}^0(C^+,\Li)^{P(\lambda)}$ onto ${\rm H}^0(C,\Li_{|C})^{G^\lambda}$.
Moreover, for any $\sigma\in{\rm H}^0(C^+,\Li)^{P(\lambda)}$, we have:
$$
\{x\in C^+\::\:\sigma(x)=0\}=\Lim^{-1}(\{x\in C\::\:\sigma(x)=0\}).
$$
\end{enumerate}
\end{prop}

\begin{proof}
Since $\Lim$ is $P(\lambda)$-equivariant, for any $\Mi\in\Pic^{G^\lambda}(C)$, 
$p_\lambda^*(\Mi)$ is $P(\lambda)$-linearized.
Since $p_\lambda$ is a vector bundle,  $p_\lambda^*(\Li_{|C})$ and $\Li$ are 
isomorphic as line bundles without linearization.
But, $\Chi(P(\lambda))\simeq\Chi(G^\lambda)$, so the 
$P(\lambda)$-linearizations must coincide; and
$p_\lambda^*(\Li_{|C})$ and $\Li$ are isomorphic as $P(\lambda)$-linearized
line bundles.
Assertion~\ref{ass1:picC+} follows.

Assertion~\ref{ass2:picC+} is a direct application of Lemma~\ref{lem:LGx}.

Let us fix $\Li\in\Pic^{P(\lambda)}(C^+)$ and denote by $p\::\:\Li\longto C^+$ the projection.
We assume that $\mu^\Li(C,\lambda)=0$.
Let $\sigma\in{\rm H}^0(C^+,\Li)^{P(\lambda)}$.
We just proved that 
$$
\Li\simeq p_\lambda^*(\Li_{|C})=
\{(x,l)\in C^+\times \Li_{|C}\::\: \Lim(x)=p(l)\}.
$$
Let $p_2$ denote the projection of $ p_\lambda^*(\Li_{|C})$ onto $\Li_{|C}$.

For all $x\in C^+$ and $t\in\kk^*$, we have:
$$
\begin{array}{lll}
\sigma(\lambda(t).x)&=\bigg(\lambda(t).x,\,p_2(\sigma(\lambda(t).x))\bigg)\\
                    &=\lambda(t).\bigg(x,\,p_2(\sigma(x))\bigg)&{\rm \ \ since\ }\sigma{\rm\ is\ invariant,}\\
                    &=\bigg(\lambda(t).x,\,p_2(\sigma(x))\bigg)&{\rm \ \ since\ }\mu^\Li(C,\lambda)=0.
\end{array}                    
$$                    
We deduce that for all $x\in C^+$, $\sigma(x)=(x,\sigma(\Lim(x)))$. 
Assertion~\ref{ass3:picC+} follows.
\end{proof}

\section{A theorem by Luna}
\label{sec:Luna}

\begin{nota}
  If $H$ is a closed subgroup of an affine algebraic group $G$, $N_G(H)$ denote the normalizer of $H$ in $G$.
\end{nota}

We will use the following interpretation of a result of Luna:

\begin{prop}
\label{prop:XssCss}
Let $\Li$ be  an ample $G$-linearized line bundle on an irreducible projective $G$-variety $X$.
Let $H$ be a reductive subgroup of $G$.
Let $C$ be an irreducible component of $X^H$.
Then, the reductive groups $(G^H)^\circ$ and $N_G(H)^\circ$ act on $C$.

Let $x$ be a point in $C$.
Then, the following are equivalent:
\begin{enumerate}
        \item $x$ is semistable for $\Li$.
        \item $x$ is semistable for the action of $(G^H)^\circ$ on $C$ and the restriction of $\Li$.
        \item $x$ is semistable for the action of $N_G(H)^\circ$ on $C$ and the restriction of $\Li$.
\end{enumerate}
\end{prop}

\begin{proof}
\cite[Lemma~1.1]{LunaRichardson:chevalley} shows that $(G^H)^\circ$ 
and $N_G(H)^\circ$ are reductive.
Changing $\Li$ by a positive power if necessary, one may assume that $X$ in contained in $\PP(V)$ 
where $V$ is a $G$-module and $\Li=\Oun_{|X}$.
Let $v\in V$ such that $[v]=x$.
Let us recall that in this case $x\in X^{\rm us}(\Li)$ if and only if $\overline{G.v}$ contains $0$.

Let $\chi$ be the character of $H$ such that $hv=\chi(h)v$ for all $h\in H$.

If $\chi$ is of infinite order, so is its restriction to the connected center $Z$ of $H$.
Then, $Z.v=\kk^*v$ and $0\in\overline{(G^H)^\circ.v}$.
In this case, $x$ belongs to no semistable set of  the proposition.

Let us now assume that $\chi$ is of finite order.
Changing $\Li$ by a positive power if necessary, one may assume that $\chi$ is trivial, that is $H$ fixes $v$.
In this case, \cite[Corollary~2 and Remark~1]{Luna:adh} shows that 
$$
0\in\overline{G.v}\iff
0\in\overline{N_G(H)^\circ.v}\iff
0\in\overline{(G^H)^\circ.v}.
$$
The proposition follows.
\end{proof}\\

\section{First descriptions of the $G$-cones}
\label{sec:Cov}

\subsection{Definitions}

Let us recall from the introduction that $\Lambda$ is a freely finitely generated subgroup of $\Pic^G(X)$ 
and $\Lambda_\QQ$ is the $\QQ$-vector space containing $\Lambda$ as a lattice.
Since $X^{\rm ss}(\Li)=X^{\rm ss}(\Li^{\otimes n})$, for any $G$-linearized line bundle and
any positive integer $n$, we can define $X^{\rm ss}(\Li)$ for any $\Li\in\Lambda_\QQ$.
The central object of this article is the following {\it total $G$-cone}:
$$
\tc^G_\Lambda(X)=\{\Li\in\Lambda_\QQ\;:\;  X^{\rm ss}(\Li) {\rm\ is\ not\ empty}\}.
$$
Since the tensor product of two non zero $G$-invariant sections is a non zero $G$-invariant section, 
$\tc^G_\Lambda(X)$ is a convex cone.

Consider the convex cones   $\Lambda^{+}_\QQ$ and $\Lambda^{++}_\QQ$ generated respectively by the semiample 
and ample elements of $\Lambda$. 
For all $\Li\in\Lambda^+_\QQ$ (resp. $\Lambda^{++}_\QQ$), there exists a positive integer $n$ such that
$\Li^{\otimes n}$ is a semiample (resp. ample) $G$-linearized line bundle on $X$ in $\Lambda$.
So, any set of semistable points associated to a point in $\Lambda^{+}_\QQ$ (resp. $\Lambda^{++}_\QQ$) 
is in fact a set  of semistable point associated to  a semiample (resp. ample) $G$-linearized line bundle.
We consider the following {\it semiample and ample $G$-cones}:
$$
\sac^G_\Lambda(X)=\tc^G_\Lambda(X)\cap\Lambda^{+}_\QQ {\rm\ \ and\ \ }   \ac^G_\Lambda(X)=\tc^G_\Lambda(X)\cap\Lambda^{++}_\QQ.
$$
By~\cite{DH} (see also \cite{GeomDedic}), $\ac^G_\Lambda(X)$ is a closed convex rational polyhedral cone in $\Lambda^{++}_\QQ$.

\subsection{Well covering pairs}
\label{sec:wellcov}

Here comes a central definition in this work:\\

\begin{defin}
Let $\lambda$ be a one parameter subgroup of $G$ and $C$ be an irreducible 
component of $X^\lambda$.
Set $C^+:=\{x\in X\ |\ \lim_{t\to 0}\lambda(t)x\in C\}$.
Consider the following $G$-equivariant map
$$
\begin{array}{cccc}
\eta\ :&G\times_{P(\lambda)}C^+&\longto& X\\
       &[g:x]&\longmapsto&g.x.
\end{array}
$$ 
The pair $(C,\lambda)$ is said to be {\it covering} (resp. {\it dominant}) if $\eta$ is birational
(resp. dominant).
It is said to be {\it well covering} if $\eta$ induces an isomorphism from 
$G\times_{P(\lambda)}\Om$ onto an open subset of $X$ 
for an open subset $\Om$ of $C^+$ intersecting $C$.
\end{defin}

Let us recall that $\mu^\bullet(C,\lambda)$ denote the common value of the 
$\mu^\bullet(x,\lambda)$, for $x\in C^+$.
The first relation between covering pairs and the $G$-cones is the following

\begin{lemma}
\label{lem:muCcov}
Let $(C,\lambda)$ be a dominant pair and $\Li\in\tc^G_\Lambda(X)$. 
Then, $\mu^\Li(C,\lambda)\leq 0$.
\end{lemma}

\begin{proof}
  Since $X^{\rm ss}(\Li)$ is a $G$-stable non empty open subset of $X$ and $(C,\lambda)$ is covering,
there exists a point $x\in C^+$ semistable for $\Li$.
Since $x\in C^+$, $\mu^\Li(x,\lambda)=\mu^\Li(C,\lambda)$.
But, since $x\in X^{\rm ss}(\Li)$ Lemma~\ref{lem:HMsemiample} shows that  $\mu^\Li(x,\lambda)\leq 0$.
\end{proof}\\

Proposition~\ref{prop:Ki} allows us to give a first description of the 
cone $\CGLX$:

\begin{prop}
\label{prop:desc1}
We assume that $X$ is normal.
Let $T$ be a maximal torus of $G$ and $B$ be a Borel subgroup containing $T$. 

Then, the cone $\ac^G_\Lambda(X)$ is the set of the $\Li\in\Lambda_\QQ^{++}$ such that for all well covering pair
$(C,\lambda)$ with a dominant one parameter subgroup $\lambda$ of $T$ we have
$\mu^\Li(C,\lambda)\leq 0$.
\end{prop}

\begin{proof}
Lemma~\ref{lem:muCcov} shows that $\ac^G(X)$ is contained is the part of $\Lambda_\QQ^{++}$ defined by the inequalities
$\mu^\Li(C,\lambda)\leq 0$ of the proposition.

Conversely, let $\Li\in\Lambda_\QQ^{++}$ such that $X^{\rm ss}(\Li)$ is empty.
Consider the open stratum $S^\Li_{d,\langle\tau \rangle }$ in $X$.
Let  $\lambda$ be a dominant one parameter subgroup of $T$ in the class $\langle\tau \rangle $.
Since $S^\Li_{d,\langle\tau \rangle }$ is open in $X$, it is normal and irreducible.
Now, Proposition~\ref{prop:Ki} implies that $Z^\Li_{d,\lambda}$ is open in an irreducible 
component $C$ of $X^\lambda$.
It shows also that  $S^\Li_{d,\lambda}$ is open in $C^+$ and intersects $C$.
Finally, the last assertion of Proposition~\ref{prop:Ki} shows that $(C,\lambda)$ is well 
covering.
Moreover, $\mu^\Li(C,\lambda)=d>0$.
\end{proof}\\

\section{Faces of the $G$-ample cone and well covering pairs}
\label{sec:gen}

\begin{defin}
Let $\varphi$ be a linear form on $\Lambda_\QQ$ which is non negative on $\CGLX$.
The set of $\Li\in\CGLX$ such that $\varphi(\Li)=0$ is  called 
{\it a face of $\CGLX$}. Note that a face can be empty or equal to $\CGLX$. 
A face different from $\CGLX$ is said to be strict.
\end{defin}

Lemma~\ref{lem:muCcov} implies that for any dominant pair $(C,\lambda)$, by intersecting $\CGLX$
with the hyperplane with equation $\mu^\bullet(C,\lambda)=0$, one obtain a face of $\CGLX$ 
(eventually empty): this face is said to be {\it associated to $(C,\lambda)$}. 
The proof of Proposition~\ref{prop:desc1} and Lemma~\ref{lem:Xcircmono} imply easily (see the 
proof of Theorem~\ref{th:facecov} below) that for any codimension one 
face $\Face$ of $\CGLX$, there exist a well covering pair $(C,\lambda)$ whose associated face 
is $\Face$. The aim of this section is to prove  that the relations between faces and well 
covering pairs are much more deeper.

\subsection{Dominant pairs and quotient varieties}

The following theorem is a description of the quotient variety associated to a point in $\CGLX$ which 
belongs to  a face associated to a dominant pair.

\begin{theo}
\label{th:CsurGlambda}
Let $(C,\lambda)$ be a dominant pair and $\Li$ be an ample $G$-linearized line bundle on $X$ such 
that $\mu^\Li(C,\lambda)=0$.
Consider the action of $G^\lambda$ on $C$ and the associated set of semistable points 
$C^{\rm ss}(\Li,G^\lambda)$. Then, $C^{\rm ss}(\Li,G^\lambda)=X^{\rm ss}(\Li,G)\cap G$. 
Consider the morphism $\theta$ which makes the following diagram commutative:
\begin{diagram}
  C^{\rm ss}(\Li_{|C},G^\lambda)&\rInto&X^{\rm ss}(\Li)\\
\dTo^{\pi_C}&&\dTo_\pi\\
C^{\rm ss}(\Li_{|C},G^\lambda)\quot G^\lambda&\rTo^\theta&X^{\rm ss}(\Li)\quot G.
\end{diagram}

Then, $\theta$ is finite and surjective.
\end{theo}

\begin{proof}
The first assertion in a direct consequence of Proposition~\ref{prop:XssCss}.
Since $C^{\rm ss}(\Li,G^\lambda)\quot G^\lambda$ is projective, to prove the second one, 
it is sufficient to prove that $\theta$ is dominant and its fibers are finite.

  Since $(C,\lambda)$ is dominant, $C^+$ must intersect $X^{\rm ss}(\Li)$.
Let $x\in C^+\cap X^{\rm ss}(\Li)$. Set $z=\lim_{t\to 0}\lambda(t)x$.
By Assertion~\ref{assGeomDedic} of Lemma~\ref{lem:HMsemiample}, $z$ is semistable for $\Li$. 
So, $C$ intersects $X^{\rm ss}(\Li)$. Moreover, we just proved that 
$\pi(C^+\cap X^{\rm ss}(\Li))= \pi(C^{\rm ss}(\Li,G^\lambda))$.
It follows that $\theta$ is dominant.

Let $\xi\in X^{\rm ss}(\Li)\quot G$.
Let $\Orb_\xi$ be the unique closed $G$-orbit in $\pi^{-1}(\xi)$.
The points in the fiber $\theta^{-1}(\xi)$ correspond bijectively to the closed $G^\lambda$-orbits in 
$\pi^{-1}(\xi)\cap C$. But, \cite[Corollary~2 and Remark~1]{Luna:adh} implies that these orbits are contained
in  $\Orb_\xi$. We conclude that $\theta^{-1}(\xi)$ is finite, 
by using a result of \cite{Rich} which implies that $\Orb_\xi^\lambda$ contains only 
finitely  many $G^\lambda$-orbits.
\end{proof}

\subsection{From faces  to well covering pairs}

\begin{nota}
  If $\lambda$ is a one parameter subgroup of $G$, $\Im\lambda$ will 
denote its image. 
If $\Face$ is  a convex cone  in a vector space, 
$\langle\Face\rangle$ will denote the subspace spanned by $\Face$.
The dimension of $\langle\Face\rangle$ will be called the dimension of $\Face$.
\end{nota}

Let $\Li$ be an ample line bundle in $\Lambda$ without semistable point. 
Let $d$ be the positive real number and $\lambda$ be a one parameter subgroup 
of $G$ such that $X^\circ(\Li)=S^\Li_{d,\langle\lambda\rangle}$. 

\begin{lemma}
  \label{lem:Clunique}
With above notation, $Z^\Li_{d,\lambda}$ (see Section~\ref{sec:strat}) is irreducible;
and the closure $C$ of $Z^\Li_{d,\lambda}$ is an irreducible component of $X^\lambda$.

The pair $(C,\lambda)$ is covering. It is well covering, if, in addition, $X$ is normal.
Moreover, the conjugacy class of the pair $(C,\lambda)$ only depends on $\Li$.

\end{lemma}

\begin{proof}
Since $S^\Li_{d,\langle\lambda\rangle}$ is open in $X$, it is irreducible.
But, $P(\lambda)$ is connected, so the last assertion of Proposition~\ref{prop:Ki} implies
that $S^\Li_{d,\lambda}$ is irreducible. It follows that $Z^\Li_{d,\lambda}$ is irreducible.
Since $Z^\Li_{d,\lambda}$ is open in $X^\lambda$ its closure is an irreducible component 
of $X^\lambda$.

The fact that $(C,\lambda)$ is covering (resp. well covering, if $X$ is normal) is a direct
consequence of Proposition~\ref{prop:Ki}.

 The last assertion is obvious since $\lambda$ is unique up to conjugacy.
\end{proof}\\

Lemma~\ref{lem:Clunique} implies that the linear 
form $\mu^\bullet(C,\lambda)$ on $\Lambda_\QQ$ only depends on $\Li$. We set:
$$
\Hcal(\Li)=\{\Li\in\Lambda_\QQ\,:\,\mu^\Li(C,\lambda)=0\}{\rm\ \ and\ \ }
\Hcal(\Li)^{>0}=\{\Li\in\Lambda_\QQ\,:\,\mu^\Li(C,\lambda)>0\}.
$$
By Proposition~\ref{prop:Ki}, the pair $(C,\lambda)$ is dominant. 
So, Lemma~\ref{lem:muCcov} implies that $\Hcal\cap \CGLX$ is a face of $\CGLX$:
this face is denoted by $\Face(\Li)$ and called the {\it face viewed from $\Li$}.

The aim of this section is to prove 

\begin{theo}
\label{th:facecov}
Any strict face $\Face$ of $\CGLX$ is viewed from some ample point in $\Lambda$ with no semistable point.
\end{theo}

We start by proving Lemmas~\ref{lem:FLp}, \ref{lem:FLconst} and \ref{lem:FLsemicont} 
about the function  $\Li\mapsto\Face(\Li)$. 
Consider the subgroup $H^{\lambda}$ of $G^{\lambda}$ defined in 
Section~\ref{sec:Hlambda}. Consider the morphism of restriction 
$p\,:\,\Lambda\longto\Pic^{H^{\lambda}}(C)$.

\begin{lemma}
\label{lem:FLp}
With above notation, $\Face(\Li)$ equals the intersection 
$\Lambda_\QQ^{++}\cap\Hcal(\lambda)\cap p^{-1}(\ac^{H^\lambda}(C))$.
\end{lemma}

\begin{proof}
Let $\Mi$ belong to the intersection of the lemma. 
There exists $x\in C$ which is semistable for the action 
of $H^{\lambda}$ and $\Li$. Since $\Mi\in \Hcal(\Li)$,  $\lambda$ acts trivially on $\Mi_{|C}$; 
so $x$ is semistable for $G^{\lambda}$.
Proposition~\ref{prop:XssCss} shows that $x$ is semistable for $G$ and $\Mi$. 
So, $\Mi\in\ac_{\Lambda}^G(X)$.

Conversely, let $\Mi\in\Face(\Li)$. 
Since $(C,\lambda)$ is covering, $X^{\rm ss}(\Mi)$ intersects $C^+$. 
But $\mu^\Mi(C,\lambda)=0$; so, Assertion~\ref{assGeomDedic} of Lemma~\ref{lem:HMsemiample} 
shows that $C$ intersects $X^{\rm ss}(\Li)$.
So, $p(\Mi)\in\ac_{p(\Lambda)}^{H^\lambda}(C)$.
\end{proof}

\begin{lemma}
\label{lem:FLconst}
Let $\Li'\in\Hcal(\lambda)^{>0}$. We assume that $p(\Li')\in \ac^{H^\lambda}(C)$. 
Then, the face of $\CGLX$ viewed from $\Li'$ equals those viewed from $\Li$. 
More precisely, the pair $(C,\lambda)$ satisfies Lemma~\ref{lem:Clunique} for $\Li'$.
\end{lemma}

\begin{proof}
  Let $x\in C^+$ such that $\lim_{t\to 0}\lambda(t)x\in C^{\rm ss}(p(\Li'),H^\lambda)$.
Note that the set of such $x$'s is open is $C^+$.
By Theorem~\ref{th:NessKi}, $\lambda$ is adapted to $x$ and $\Li'$.
Since $(C,\lambda)$ is covering, $S^{\Li'}_{\overline{\mu}^{\Li'}(x,\lambda),\langle\lambda\rangle}=
X^\circ(\Li')$. The lemma follows.  
\end{proof}

\begin{lemma}
\label{lem:FLsemicont}
  There exists an open neighborhood $U$ of $\Li$ such that for any $\Li'$ in $U$, $\Li'$ has
no semistable point and the face of $\CGLX$ viewed from $\Li'$ is contained in those viewed 
from $\Li$.

If in addition $p(\Li')$ has no semistable point, we have $\Face(\Li')$ is the set of 
$\Mi\in\Hcal(\lambda)$ such that $p(\Mi)$ belongs to the face of $\ac^{H^\lambda}_{p(\Lambda)}(C)$ 
viewed from $p(\Li')$.
\end{lemma}

\begin{proof}
  Let $T$ be a maximal torus of $G$ containing the image of $\lambda$.
By Lemma~\ref{lem:convex}, there exists an open neighborhood $U$ of $\Li$ such that for any 
$\Li'$ in $U$ we have:
\begin{center}
  For any $x\in C^+\cap X^\circ(\Li)$, the face of $\P_T^{\Li'}(x)$ viewed from $0$
  is contained in the hyperplane $\langle \lambda,\cdot\rangle=\mu^{\Li'}(x,\lambda)$.
\end{center}
Let us fix $\Li'\in U$ and $x\in C^+\cap X^\circ(\Li)$ such that 
$z=\lim_{t\to 0}\lambda(t).x$ belongs to $C^\circ(\Li',H^\lambda)$. 
By \cite[Lemma~12.19]{Ki}, there exists a one parameter subgroup adapted to $z$ and $\Li'$ 
which commutes with $\lambda$. So, there exists  $h_0\in H^\lambda$ and a one parameter subgroup $\zeta$ 
of $T$ which is adapted to $h_0z$ and $\Li'$. 
Since $h_0x\in C^+\cap X^\circ(\Li)$, the orthogonal projection of $0$ on $\P^{\Li'}_T(h_0x)$ belongs
to $\P_T^{\Li'}(h_0z)$. In particular, 
$\lim_{t\to 0}\zeta(t)h_0x=\lim_{t\to 0}\zeta(t)h_0z=:z'$. 
Let $H^\zeta$ the subgroup of $G^\zeta$ defined in Section~\ref{sec:Hlambda}.
Since $\zeta$ is adapted to $h_0z$ and $\Li'$,  Theorem~\ref{th:NessKi} implies that $z'$ is
semistable for the action $H^\zeta$ and $\Li'$. So, by Theorem~\ref{th:NessKi} again, $\zeta$ 
is adapted for $h_0x$ and $\Li'$.

One easily checks that $C^\circ(\Li',H^\lambda)=C^\circ(\Li',G^\lambda)$.
Let $d$ be such that $ C^\circ(\Li',G^\lambda)=S^{\Li'}_{d',\langle\zeta\rangle}$.
So, for any $z'\in C^\circ(\Li',H^\lambda)$, there exists $h'\in H^\lambda$ such that 
$\zeta$ is adapted to $h'z'$ and $\Li'$. So, we just proved that for any $x$ in a nonempty
open subset of $C^+$, $\zeta$ is adapted to $h.x$ for some $h\in H$. 
Since $G.C^+$ contains an open subset of $X$, we can deduce that 
$X^\circ(\Li')=S^{\Li'}_{d',\langle\zeta\rangle}$. 
Note that the equality $\lim_{t\to 0}\zeta(t)h_0x=\lim_{t\to 0}\zeta(t)h_0z$ implies that 
$Z^{\Li'}_{d',\zeta}$ is contained in $C$. Let  $C_\zeta$ denote the closure of  $Z^{\Li'}_{d',\zeta}$.

Since $\P_T^{\Li'}(h_0z)$ is contained in the hyperplane 
$\langle\lambda,\cdot\rangle=\mu^{\Li'}(z,\lambda)>0$, there exists a one parameter subgroup $\lambda_1$
of $T\cap H^\lambda$ such that $\zeta=\overline{\lambda+\lambda_1}$.
Since $z$ is semistable for the action of $H^\lambda$ and $\Li$, $0$ belongs to 
$\P_{T\cap H^\lambda}^{\Li}(h_0z)$. We deduce that $\lambda_1$ tends to 0 when $\Li'$ tends to $\Li$.

Since $C_\zeta\subset C$, there exists a positive constant $c$ such that:
$$
\mu^\bullet(C_\zeta,\zeta)=c.\mu^\bullet(C_\zeta,\lambda_1+\lambda)=
c.(\mu^\bullet(C_\zeta,\lambda_1)+\mu^\bullet(C,\lambda)).
$$
In particular, $\mu^\bullet(C_\zeta,\zeta)$ tends to $\mu^\bullet(C,\lambda)$  
when $\Li'$ tends to $\Li$.
Up to changing $U$, we may assume that $\Face(\Li')\subset \Face(\Li)$ for any $\Li'\in U$.

Now, Lemma~\ref{lem:FLp} implies that for any $\Li'\in U$, $\Face(\Li')$ is the set of $\Mi\in\Lambda^+$
such that $\Mi\in\Hcal(\lambda)$ and $p(\Mi)\in\Face(p(\Li'))$. The lemma is proved. 
\end{proof}\\

We are now ready to prove the theorem.\\

\begin{proof}[of Theorem~\ref{th:facecov}]
We will prove the following assertion, by induction on the codimension $\codim_\Lambda(\Face)$ 
of $\Face$ in $\Lambda$:

{\it Let $X$ be a not necessarily normal variety. 
Let $G$, $\Lambda$ and $\Face$ as in the theorem. 
Let $U$ be an open subset of $\Lambda_\QQ$ intersecting $\Face$.
Then, there exists $\Li\in U$ such that $\Face(\Li)=\Face$.}\\

Let $\Mi$ be a point in the relative interior of $\Face\cap U$.
By Lemma~\ref{lem:Xcircmono}, there exists $\Li\in U$ without semistable point and such that
$X^{\rm \circ}(\Li)\subset X^{\rm ss}(\Mi)$.
Let $(C,\lambda)$ associated to $\Li$ as in Lemma~\ref{lem:Clunique}.
Since $C$ intersects $X^{\rm \circ}(\Li)$ it intersect $X^{\rm ss}(\Mi)$; so, 
$\mu^\Mi(C,\lambda)=0$ and $\Mi\in\Face(\Li)$.
Since $\Mi$ belongs to the relative interior of $\Face$, we deduce that $\Face\subset\Face(\Li)$. 

If $\codim_\Lambda(\Face)=1$ then   $\Face=\Face(\Li)$, and the theorem is proved.
Let us now assume that $\codim_\Lambda(\Face)\geq 2$. 

By Lemma~\ref{lem:FLp}, $\Face\cap U$ is a face
 of $U\cap\Hcal(\Li)\cap p^{-1}(\ac^{H^\lambda}_{p(\Lambda)}(C))$. 
So, there exists a face $\tilde\Face_1$ of $p^{-1}(\ac^{H^\lambda}_{p(\Lambda)}(C))$
such that $\Face=\Lambda^{++}_\QQ\cap\Hcal(\Li)\cap\tilde\Face_1$.
Let us assume that $\tilde\Face_1$ is of maximal dimension among such faces. 
Since $\Hcal(\Li)^{>0}\cap p^{-1}(\ac^{H^\lambda}_{p(\Lambda)}(C))$ is non empty (it contains $\Li$ !),
we may assume that   $\tilde\Face_1$ intersects $\Hcal(\Li)^{>0}$.
Let $\Face_1$ denote the unique face of $\ac^{H^\lambda}_{p(\Lambda)}(C)$ such that 
$\tilde\Face_1=p^{-1}(\Face_1)$.

Since $\tilde\Face_1$ intersects $\Hcal(\Li)^{>0}$, Lemma~\ref{lem:FLconst} allows to move $\Li$ 
on $\tilde\Face_1$ without changing $\Face(\Li)$ nor $(C,\lambda)$. From now on, we assume that 
$\Li\in\tilde\Face_1$.

Up to restricting  $U$ if necessary, we may assume it  satisfies Lemma~\ref{lem:FLsemicont} and 
it is contained in $\Hcal(\Li)^{>0}$.
Since $\Face\subset \tilde\Face_1\cap \Hcal(\Li)$ and  $\tilde\Face_1$ intersects $\Hcal(\Li)^{>0}$,
we have: $\codim_\Lambda(\tilde\Face_1)<\codim_\Lambda\Face$. 
But, $\codim_\Lambda(\tilde\Face_1)=\codim_{p(\Lambda)}(\Face_1)$.
Moreover, $p$ is linear and so open. We apply the induction to the action of $H^\lambda$ on $C$,
the face $\Face_1$ and the open subset $p(U)$ of $p(\Lambda)_\QQ$: there exists $\Li'$ in $U$ such that 
$\Face(p(\Li'))=\Face_1$.
Lemma~\ref{lem:FLsemicont} implies now that $\Face(\Li')=\Face$.
\end{proof}\\

Theorem~\ref{th:facecov} can be restated in terms of well covering pairs:

\begin{coro}\label{cor:facecov}
  Let $\Face$ be a face of $\CGLX$ of codimension $r$ in $\Lambda$.
Then, there exists a $r$-dimensional torus $S$ in $G$, a one parameter subgroup $\lambda$
of $S$ and an irreducible component $C$ of $X^S$ such that:
\begin{enumerate}
\item $C$ is an irreducible component of $X^\lambda$;
\item the pair $(C,\lambda)$ is covering (respectively, well covering if $X$ is normal);
\item $\Face$ is the face of $\CGLX$ associated to $(C,\lambda)$.
\end{enumerate}
\end{coro}

\begin{proof}
 Let $\Li$ be an ample line bundle on $X$ without semistable point and such that $\Face$ 
is viewed from $\Li$. Let $d$ be the positive number and $\lambda$ be a one parameter 
subgroup of $G$ such that $S^\Li_{d,\langle\lambda\rangle}$ is open in $X$.
Let $C$ be the closure of $Z^\Li_{d,\lambda}$.

By Lemma~\ref{lem:Clunique}, it remains to prove that there exits a torus $S$ of dimension $r$
containing the image of $\lambda$ and acting trivially on $C$.
The proof of the existence of such a $S$ can easily be integrated in the induction of the proof 
of Theorem~\ref{th:facecov}.
\end{proof}

\subsection{From well covering pairs to faces of total $G$-cone}

In Corollary~\ref{cor:facecov}, starting with a face of $\ac^G(X)$ we have constructed a well covering pair.
Conversely, in the following theorem, we start with a well covering pair and study the associated face of $\tc^G(X)$.

\begin{theo}
\label{th:wc2face}
Let $X$ be a projective $G$-variety.
We assume the rank of  $\Pic^G(X)$ is finite and consider $\tc^G(X)$.
Let $(C,\lambda)$ be a well covering pair.
Consider the  linear map $\rho$ induced by the restriction:
$$
\rho\ :\ \Pic^G(X)_\QQ\longto\Pic^{G^\lambda}(C)_\QQ.
$$
Then, the subspace of $\Pic^G(X)_\QQ$ spanned by the $\Li\in\tc^G(X)$ such that
$\mu^\Li(C,\lambda)=0$ is the pullback by $\rho$ of the intersection of the image of $\rho$ 
and  $\tc^{G^\lambda}(C)$.
\end{theo}

\begin{proof}
  Let $\Li\in\tc^G(X)$ such that $\mu^\Li(C,\lambda)=0$.
Since $X^{\rm ss}(\Li)$ is open and $G$-stable, and $(C,\lambda)$ is covering,
there exists $x\in C^+$ semistable for $\Li$.
Set $z=\lim_{t\to 0}\lambda(t)x$. 
Since $\mu^\Li(x,\lambda)=0$,  Lemma~\ref{lem:HMsemiample} shows that $z$ is 
semistable for $\Li$. This implies that the restriction of $\Li$ belongs to
$\tc^{G^\lambda}(C)$.
This proves that the first subspace of the theorem is contained in the 
second one.\\

We denote by $F$ the  pullback by $\rho$ 
of the subspace spanned by $\tc^{G^\lambda}(C)$.
Let $\Li_1,\,\dots,\,\Li_n\in\Pic^G(X)$ which span $F$ and  whose the restrictions 
to $C$ belong to $\tc^{G^\lambda}(C)$.
Denote by $\Mi_i$ the restriction of $\Li_i$ to $C$. 
For each $i$, let us fix a non zero regular $G^\lambda$-invariant section 
$\sigma_i$ of $\Mi_i$.
  
Consider the $G$-linearized line bundle $G\times_{P(\lambda)}\Lim^*(\Mi_i)$ on 
$G\times_{P(\lambda)}C^+$ with notation of Section~\ref{sec:BB}.
Since $\eta^*(\Li_i)$ and $G\times_{P(\lambda)}\Lim^*(\Mi_i)$ have the same restriction to 
$C$, Lemma~\ref{lem:picGPY} and Proposition~\ref{prop:picC+} show that 
$\eta^*(\Li_i)=G\times_{P(\lambda)}\Lim^*(\Mi_i)$.
Moreover, since $\mu^{\Mi_i}(C,\lambda)=0$, Proposition~\ref{prop:picC+} shows that 
$\sigma_i$ admits a unique $P(\lambda)$-invariant extension to a section $\sigma_i'$
of $\Lim^*(\Mi_i)$. 
On the other hand, Lemma~\ref{lem:picGPY} shows that $\sigma_i'$ admits a unique
$G$-invariant extension $\tilde{\sigma}_i$ from $C^+$ to $G\times_{P(\lambda)}C^+$.
So, we obtain the following commutative diagram:

\begin{diagram}
\Li_i&\lTo&\eta^*(\Li_i)=G\times_{P(\lambda)}\Lim^*(\Mi_i)&\lTo&\Lim^*(\Mi_i)&\rTo&\Mi_i\\
\dTo&   &\dTo\uTo_{\tilde{\sigma}_i}&            &\dTo\uTo_{\sigma_i'}&&\dTo\uTo_{\sigma_i}\\
X&\lTo^\eta&G\times_{P(\lambda)}C^+&\lTo&C^+&\rTo^\Lim&C
\end{diagram}

Since $\eta$ is birational, $\tilde{\sigma}_i$ descends to a rational $G$-invariant 
section $\tau_i$ of $\Li_i$.
Let $X^\circ$ be a $G$-stable open subset of $X$ such that $\eta$ induces an isomorphism from
$\eta^{-1}(X^\circ)$ onto $X^\circ$. Since $(C,\lambda)$ is well covering, we may (and shall)
assume that $X^\circ$ intersects $C$.
Let $E_j$ be the irreducible components of codimension one of $X-X^\circ$.
For any $j$ we denote by $a_j$ the maximum of $0$ and the $-\nu_{E_j}(\tau_i)$'s with 
$i=1,\,\dots,\,n$.
Consider the line bundle $\Li_0=\O(\sum a_jE_j)$ on $X$.
Since the $E_j$'s are stable by the action of $G$, $\Li_0$ is canonically $G$-linearized.
By construction, the $\tau_i$'s induce $G$-invariant regular sections $\tau_i'$ of 
$\Li_i':=\Li_i\otimes\Li_0$.
Moreover, since no $E_j$ contains $C$, the restriction of $\tau_i'$ to $C$ is non zero.
In particular, the $\Li_i'$'s belong to $\tc^G(X)$ and their restrictions belong to $\tc^{G^\lambda}(C)$.
Moreover, replacing $\Li_0$ by $\Li_0^{\otimes 2}$ if necessary, we may (and shall) assume that the 
 $\Li_i'$'s span $F$.
This ends the proof of the theorem.
\end{proof}

\subsection{If $X=G/B\times Y$\dots}

In this section, we will explain how Theorems~\ref{th:CsurGlambda} and \ref{th:facecov} can be 
improved when $X=G/B \times Y$ for a $G$-variety $Y$.

\paragraph {\bf Complement to Theorem~\ref{th:CsurGlambda}.}

\begin{theo}
\label{th:CxcGB}
We use notation of Theorem~\ref{th:CsurGlambda}, assuming in addition that
$X=G/B\times Y$ with a normal projective $G$-variety $Y$.
We also assume that $X^{\rm ss}(\Li)$ is not empty.

Then, $\theta$ is an isomorphism and $(C,\lambda)$ is a well covering pair.
\end{theo}

\begin{proof}
Since $X$ is normal, so is $X^{\rm ss}(\Li)\quot G$. But, Theorem~\ref{th:CsurGlambda} shows
that $\theta$ is finite; it is sufficient to prove that it is birational. Since the base field
 has characteristic zero, it is sufficient to prove that $\theta$ is bijective to obtain the 
first assertion.

Let $\xi\in X^{\rm ss}(\Li)\quot G$ and $\Orb$ denote the only closed $G$-orbit 
in $\pi^{-1}(\xi)$.
As noticed during the proof of  Theorem~\ref{th:CsurGlambda}, $\theta^{-1}(\xi)$
correspond bijectively with the set of $G^\lambda$-orbits in $\Orb\cap C$.

Consider the first projection $p_1\,:\,X\longto G/B$ and fix $x\in\Orb\cap C$.
Let $B$ denote the stabilizer in $G$ of $p_1(x)$.
Since $\Orb$ is closed in $\pi^{-1}(\xi)$, it is affine and so, $G_x$ is reductive. But,
$G_x$ is contained in $B$, so $G_x$ is diagonalisable. Let $T$ be a maximal torus of 
$G$ such that $G_x\subset T\subset B$.
Consider the $G$-equivariant morphisms $q\,:\,\Orb\longto G/T,\,x\mapsto T/T$ and 
$\hat q\,:\,G/T\longto G/B,\,T/T\mapsto B/B$ unduced by these
inclusions; we have, $p_1=\hat q\circ q$.

We claim that ${\hat q}^{-1}(G^\lambda B/B)^\lambda=G^\lambda T/T$.
Since $G^\lambda$ is connected, each irreducible component of 
 ${\hat q}^{-1}(G^\lambda B/B)^\lambda$ is $G^\lambda$-stable, and so, it maps onto $G^\lambda T/T$.
In particular, it intersects ${\hat q}^{-1}(B/B)=B/T$.
But, $B/T$ is isomorphic to the Lie algebra of the unipotent radical of $B$ as a $T$-variety;
in particular, $(B/T)^\lambda$ is irreductible and so is ${\hat q}^{-1}(G^\lambda B/B)^\lambda$.
The claim now follows from \cite{Rich}.

Note that $C$ is the product of one irreducble component of $(G/B)^\lambda$ and one of 
$Y^\lambda$; this implies that  $p_1(C)=G^\lambda B/B$.
So, $\Orb\cap C$ is contained in $p_1^{-1}(G^\lambda B/B)^\lambda$ and so in 
$\Orb\cap q^{-1}(G^\lambda T/T)$.
But, since $T$ is contained in $G^\lambda$, $q^{-1}(G^\lambda T/T)=G^\lambda x$.
It follows that $\Orb\cap C=G^\lambda .x$; and so, that $\theta^{-1}(\xi)$ is reduced 
to one point. This ends the proof of the first assertion.\\

Consider $\eta\,:\,G\times_{P(\lambda)}C^+\longto X$.
We claim that for any $x\in C^+\cap X^{\rm ss}(\Li)$, $\eta^{-1}(x)$ 
is only one point.
Since $\eta$ is dominant and  the ground field has characteristic zero, the
claim implies that $\eta$ is birational.
Since $X^{\rm ss}(\Li)\cap C$ is non empty, the claim implies that $\eta$ is 
bijective over on open subset of $X$ intersecting $C$. 
Since $X$ is normal, Zariski's main Theorem (see~\cite[\S 8.12]{ega43}) implies that $(C,\lambda)$
is well covering.

Let us prove the claim.
Let $g\in G$ such that $g^{-1}x\in C^+$. We have  to prove that $g\in P(\lambda)$.
Set $x'=g^{-1}x$, $z=\lim_{t\to 0}\lambda(t)x$ and $z'=\lim_{t\to 0}\lambda(t)x'$.
Obviously $x'\in X^{\rm ss}(\Li)$, and by  Assertion~\ref{assGeomDedic} of 
Lemma~\ref{lem:HMsemiample}, $z,\,z'\in X^{\rm ss}(\Li)$ as well.  
It is also clear that $\pi(x)=\pi(x')=\pi(z)=\pi(z')=:\xi$.

Let $x_0\in C\cap\pi^{-1}(\xi)$ whose the orbit is closed in $X^{\rm ss}(\Li)$. Set $H=G_{x_0}$.
Consider $\Sigma$ (resp. $\Sigma_C$) the set of $y$ in $X$ (resp. $C$) such that $x_0$
is contained in the closure of $H.y$ (resp. $H^\lambda.y$).
By \cite{Luna:slice} (see also \cite{PV}), 
$\pi^{-1}(\xi)\simeq G\times_H\Sigma$ and $\pi_C^{-1}(\theta^{-1}(\xi))\simeq G^\lambda\times_{H^\lambda}
\Sigma_C$.
Consider the natural $G$-equivariant morphism $\gamma\,:\,G\times_H\Sigma\longto G/H$.
Since $\pi^{-1}(\xi)\cap C=\pi_C^{-1}(\theta^{-1}(\xi))$, it equals $G^\lambda.\Sigma_C$.
So, $\gamma(C\cap \pi^{-1}(\xi))=G^\lambda H/H$.
Since $\gamma$ is ``continuous'' and $G$-equivariant, we deduce that $\lim_{t\to 0}\lambda(t)\gamma(x)$ and 
 $\lim_{t\to 0}\lambda(t)\gamma(x')$ belong to $G^\lambda H/H$.
Lemma~\ref{lem:Rich+} below proves that there exists $p$ and $p'$ in $P(\lambda)$
such that $\gamma(x)=pH/H$ and  $\gamma(x')=p'H/H$. Since $gx'=x$, we have $g\in pHp'^{-1}$.
But, $H$ is a reductive subgroup (since $G.x_0$ is closed in the affine variety $\pi^{-1}(\xi)$) 
contained in a Borel subgroup of $G$ (since $X=G/B\times Y$) and containing the image of $\lambda$: 
it follows that $H\subset G^\lambda$.
Finally, $g\in P(\lambda)$.  
\end{proof}\\

It remains to prove the following lemma. In fact, it is an adaptation of the main result of \cite{Rich}:

\begin{lemma}
\label{lem:Rich+}
Let $\Orb$ be any $G$-homogeneous space. 
Let $\lambda$ be a one parameter subgroup of $G$. 
Let $C$ be an irreducible component of fixed points of $\lambda$ in $\Orb$.
Set $C^+:=\{x\in \Orb\,:\,\lim_{t\to 0}\lambda(t)x {\rm\ exists\ and\ belongs\ to\ }C\}$.

Then, $C^+$ is a $P(\lambda)$-orbit.
\end{lemma}

\begin{proof}
  Let $x\in C$. The differential of the map $g\mapsto g.x$ induces a 
surjective linear map $\phi\,:\,\lg\longto T_x\Orb$.
Since $x$ is fixed by $\lambda$, $\lambda$ acts on $T_x\Orb$; it also acts by 
the adjoint action on $\lg$, and $\phi$ is equivariant.
In particular, the restriction of $\phi$, 
$\bar\phi\,:\,\lg_{\geq 0}\longto (T_x\Orb)_{\geq 0}$ is also surjective.
But, on one hand $(T_x\Orb)_{\geq 0}=T_xC^+$ by Theorem~\ref{th:BB}  and
on the other hand $\lg_{\geq 0}$ is the Lie algebra of $P(\lambda)$. 
One can conclude that $T_xC^+=T_xP(\lambda)x$.
Since $P(\lambda)x$ is smooth, this implies that $P(\lambda)x$ is open in $C^+$.
Since $C^+$ is irreducible, it contains a unique open 
$P(\lambda)$-orbit $\Orb_0$ which contains $C$.

Note that $C^+-\Orb_0$ is $P(\lambda)$-stable and closed in $C^+$.
Since for any $y\in C^+$ $\lim_{t\to 0}\lambda(t)y\in C\subset\Orb_0$, this implies that 
$C^+-\Orb_0$ is empty.  
\end{proof}\\

\begin{remark}
Here, is an example which proves that the assumption on $X$ is useful.
Let $V$ be a vector space of dimension 2. 
Make the group $G=\kk^*\times {\rm SL}(V)$ acting on $X=\PP(\kk\oplus V\oplus V)$ by:
\begin{eqnarray}
  \label{eq:exple}
(t,g).[\tau:v_1:v_2]=[\tau,gv_1,t^2g v_2].
\end{eqnarray}

Consider (with obvious notation), the following one parameter subgroup 
$\lambda(t)=(t^{-1},\left(
  \begin{array}{cc}
    t^{-1}&0\\0&t
  \end{array}\right)$.
Note that $C=\{[1:0:0]\}$ is an irreducible component of $X^\lambda$. 
Formula~\ref{eq:exple} gives also a linearization $\Li$ of the bundle ${\mathcal O}(1)$.
Then, $\mu^\Li(C,\lambda)=0$, the map $\eta$ associated to $(C,\lambda)$ is birational, but the fiber over the only point of $C$ is $\PP^1$:
so $(C,\lambda)$ is covering but not well covering. 
\end{remark}

\paragraph {\bf Our assumption on $\Lambda$.} 
We call a subgroup $\Gamma'$ of an abelian group $\Gamma$ {\it cofinite} if $\Gamma/\Gamma'$ is finite.
The following definition is an adaptation of those of Dolgachev and Hu (see \cite{DH}).\\

\begin{defin}
The subgroup $\Lambda$ is said to be {\it abundant} if for any $x$ 
in $X$ such that $G_x$ is reductive, the image of the restriction 
$\Lambda\longto\Pic^G(G.x)$ is cofinite.
\end{defin} 

The main example of abundant subgroups comes from the case when $X=G/B\times Y$.

\begin{prop}
\label{prop:YGBabu}
Let $X=G/B\times Y$ for a  $G$-variety $Y$.
Let $\pi\ :\ X\longto G/B$ denote the projection map and 
$\pi^*\ :\ \Pic^G(G/B)\longto\Pic^G(X)$ the associated homomorphism.

Then, any subgroup $\Lambda$ of $\Pic^G(X)$ containing the image of $\pi^*$ is
abundant.
\end{prop}

\begin{proof}
Let $x=(y,gB/B)\in G/B\times Y$.
Changing $x$ by $g^{-1}x$, we assume that $g=e$.
Let $\chi\in\Chi(G_x)$. 
Note that $G_x=B_y$. Since the restriction map $\Chi(B)\longto\Chi(B_y)$ is surjective, 
there exists $\nu\in\Chi(B)$ such that $\nu_{|B_y}=\chi$.
The restriction of $\pi^*(\Li_\nu)$ to $G.x$ equals $\Li_\chi$; the proposition follows.
\end{proof}\\

\paragraph{} The next statement is a precision of Corollary~\ref{cor:facecov}.

\begin{coro}\label{cor:facecovGB}
Let $X=G/B\times Y$ for a  normal $G$-variety $Y$.
We assume that $\Lambda\subset\Pic^G(X)$ is abundant.
Let $\Face$ be a face of $\CGLX$ of codimension $r$ in $\Lambda$.
Then, there exists a $r$-dimensional torus $S$ in $G$, a one parameter subgroup $\lambda$
of $S$ and an irreducible component $C$ of $X^S$ such that:
\begin{enumerate}
\item for generic $x\in C$, we have $G_x^\circ=S$; 
\item $C$ is an irreducible component of $X^\lambda$;
\item the pair $(C,\lambda)$ is covering (respectively, well covering if $X$ is normal);
\item $\Face$ is the face of $\CGLX$ associated to $(C,\lambda)$.
\end{enumerate}
\end{coro}

\begin{proof}
   We just have to prove that if $(C,\lambda)$ and $S$ satisfy  Corollary~\ref{cor:facecov}, 
they also satisfy the first assertion.

Let $\Li$ be an ample $G$-linearized line bundle on $X$ and $x$ be a semistable point for $\Li$.
We claim that $G_x$ is diagonalizable. Consider $\pi\,:\,X^{\rm ss}(\Li)\longto X^{\rm ss}(\Li)\quot G$.
Let $y$ be a point in the closure of $G.x$ such that $G.y$ is closed in $X^{\rm ss}(\Li)$.
Since $G.y$ is closed in $\pi^{-1}(\pi(y))$, it is affine; so, $G_y$ is reductive.
But $G_y$ is contained in a Borel subgroup of $G$; so, $G_y$ is diagonalizable.
By  Luna's Slice Theorem, $\pi^{-1}(\pi(y))\simeq G\times_{G_y}\Sigma$ for an affine 
$G_y$ variety $\Sigma$. Since $x\in \pi^{-1}(\pi(y))$, we deduce that $G_x$ is conjugated to 
a subgroup of $G_y$. In particular, $G_x$ is diagonalizable.

Let us fix a point $x$ in the (finite) intersection of the set $X^{\rm ss}(\Li)\cap C$ 
for $\Li\in\Face$. 
Consider $\rho\,:\,\Lambda\longto X(G_x)$. 
Since $\Lambda$ is abundant and $G_x$ diagonalizable, the rank of $\rho$ equal the dimension of 
$G_x$. Since $\Face$ is contained in the kernel of $\rho$, the dimension of $G_x$ is less or 
equal to $r$.
Since $S\subset G_x$, it follows that $G_x^\circ=S$.
\end{proof}

\section{If $X=G/Q\times \hG/\hQ$\dots}
\label{sec:BeSj}

\subsection{Interpretations of the $G$-cones}

From now on, we assume that $G$ is a connected reductive subgroup of a connected reductive group $\hG$.
Let us fix  maximal tori $T$ (resp. $\hT$) and Borel subgroups $B$ (resp. $\hB$) of $G$ (resp. $\hG$) such that 
$T\subset B\subset \hB\supset\hT\supset T$. Let $Q$ (resp. $\hQ$) be a parabolic subgroup 
of $G$ (resp. $\hG$) containing $B$ (resp. $\hB$); let $L$ (resp. $\hL$) denote the Levi subgroup of $Q$ 
(resp. $\hQ$) containing $T$ (resp. $\hT$).

In this section, $X$ denote the variety $G/Q \times \hG/\hQ$ endowed with the diagonal action of $G$.
We will apply the results of Section~\ref{sec:gen} to $X$ with $\Lambda=\Pic^G(X)$.
The cones $\tc^G(X),\,\sac^G(X)$ and $\ac^G(X)$ will be denoted without the $\Lambda$ in subscribe.

Let us describe  $\Pic^G(X)_\QQ$.
Consider the natural action of $G\times\hG$ on $X$.
By applying the construction of Section~\ref{sec:PicGP} to the $G\times\hG$-homogeneous space $X$,
one obtains the following isomorphism $\Chi(Q)\times\Chi(\hQ)\longto \Pic^{G\times\hG}(X),
\,(\nu,\hnu)\longmapsto \Li_{(\nu,\,\hnu)}$.

\begin{lemma}
\label{lem:PicGQ}
The following short sequence is exact
$$
\begin{diagram}
0&\rTo&\Chi(\hG)_\QQ&\rTo^{\hnu\mapsto\Li_{(-\hnu_{|Q},\,\hnu_{|\hQ})}}& \Pic^{G\times\hG}(X)_\QQ&\rTo^{r^{\Delta G}_\QQ}&\Pic^G(X)_\QQ&\rTo&0\\
\end{diagram}
$$
where the second linear map is induced by the  restriction $r^{\Delta G}$ of the action of 
$G\times\hG$ to $G$ diagonally embedded in $G\times\hG$.
\end{lemma}

\begin{proof}
Let $\chi$ and $\hchi$ be characters of $G$ and $\hG$ respectively.
The trivial bundle on $X$ linearized by $(\chi,\hchi)$ belongs to
$\Pic^{G\times\hG}(X)$.
The image of this line bundle in $\Pic^G(X)$ is the trivial line bundle
linearized by the character $\chi+\hchi_{|G}$ of $G$.
In particular, any $G$-linearization of the trivial bundle belongs to the image of $r^{\Delta G}$.

Let $\Li\in\Pic^G(X)$. Let $\Li'\in\Pic(X)$ obtained from $\Li$ by forgetting the action of $G$.
By \cite{KKV}, there exists a $G\times\hG$-linearization $\Mi$ of $\Li'^{\otimes n}$ for a positive 
integer $n$.
Then $\Mi^*\otimes\Li^{\otimes n}$ is the trivial line bundle over $X$; so, it belongs to the 
image $r^{\Delta G}$.  
Finally, $\Li^{\otimes n}$ belongs to the image of the restriction.
This ends the proof of the surjectivity of $r^{\Delta G}_\QQ$.\\

Let $\Li$ in the kernel of $r^{\Delta G}$.
Since $\Li$ is trivial as a line bundle, there exists characters $\chi$ and $\hchi$ of $G$ and $\hG$ 
such that $\Li=\Li_{(\chi,\,\hchi)}$. The $G$-linearization of this last character is trivial if
and only if $\chi+\hchi_{|G}$ is trivial. This ends the proof of the lemma.
\end{proof}\\

\begin{prop}
\label{prop:interp}
\begin{enumerate}
\item $\tc^G(X)=\sac^G(X)$ is a closed convex polyhedral cone.
\item Let $(\nu,\hnu)\in\Chi(Q)_\QQ\times\Chi(\hQ)_\QQ$.
Then $r^{\Delta G}_\QQ(\Li_{(\nu,\hnu)})\in \tc^G(X)$ if and only if  
$\hnu$ and $\nu$ are dominant and for $n$ big enough $V_{n\nu}\otimes V_{n\hnu}$ 
contains non zero $G$-invariant vectors.
\item If $\ac^G(X)$ is non empty, its closure in $\Pic^G(X)_\QQ$ is $\sac^G(X)$.
\item If $Q$ and $\hQ$ are Borel subgroups of $G$ and $\hG$ then $\ac^G(X)$ is non empty.
\end{enumerate}
\end{prop}

\begin{proof}
  Since $X$ is homogeneous under the action of $G\times\hG$ and since any line bundle has a linearizable power, 
every bundle with a non zero section is semiample. Hence, $\tc^G(X)=\sac^G(X)$.
Lemma~\ref{lem:HMsemiample} and Proposition~\ref{prop:Mfini} implies that 
only finitely many non empty open subsets of $X$ are of the form  $X^{\rm ss}(\Li)$ with 
$\Li\in\Pic^G(X)$.
Let $x$ be a point in the intersection of these subsets. 
Lemma~\ref{lem:HMsemiample} shows that $\sac^G(X)$ is the set of points $\Li\in\Pic^G(X)_\QQ$ semiample
such that $\mu^\Li(x,\lambda)\leq 0$ for all one parameter subgroup $\lambda$ of $G$. 
In particular, $\sac^G(X)$ is a closed convex cone.\\

Let now $\nu$ and $\hnu$ be characters of $Q$ and $\hQ$.
If $\nu$ or $\hnu$ is not dominant then $\Li_{(\nu,\hnu)}$ has no regular section.
Else, Borel-Weil's theorem shows that ${\rm H}^0(X,\Li_{(\nu,\hnu)})$ is isomorphic 
as a $G\times\hG$-module to $V_\nu^*\otimes V_\hnu^*$. 
The second assertion of the proposition follows.

The third assertion is satisfied since the ample cone in $\Pic^G(X)$ is the interior 
of the semiample cone.\\

Let $w_0$ be the longest element of $W$.
Assume that $Q=B$ and $\hQ=\hB$ are Borel subgroups of $G$ and $\hG$.
Let $\hnu_0$ be any character of $\hB$ such that $\Li_{\hnu_0}$ is ample over $\hG/\hB$.
Let $\nu$ be any dominant weight of the $G$-module $V_{\hnu_0}^*$. 
Then, $r^{\Delta G}(\Li_{(-w_0\nu,\hnu_0)})$ belongs to $\sac^G(X)$.

Let $\nu_0$  be any character of $B$ such that $\Li_{\nu_0}$ is ample over $G/B$.
Since the restriction $\kk[\hG]\longto\kk[G]$ is surjective, Frobenius's theorem implies
that $V_{\nu_0}^*$ is contained in an irreducible $\hG$-module $V_{\hnu}^*$. 
Then, $r^{\Delta G}(\Li_{(-w_0\nu_0,\hnu)})$ belongs to $\sac^G(X)$.

Since $\sac^G(X)$ is convex, it contains $r^{\Delta G}(\Li_{(-w_0(\nu_0+\nu),\hnu+\hnu_0)})$.
But the line bundle  $\Li_{(-w_0(\nu_0+\nu),\hnu+\hnu_0)}$ is ample. The last assertion is proved.
\end{proof}\\

\subsection{Covering and well covering pairs}

\begin{nota}
Let $W$ and $\hW$ denote the Weyl groups of $G$ and $\hG$. 
If $P$ is a parabolic subgroup of $G$ containing $T$, $W_P$ denote the Weyl group of the Levi
subgroup of $P$ containing $T$. 
This group $W_P$ is canonically embedded in $W$.  
\end{nota}

In this subsection, we fix two parabolic subgroups $Q$ and $\hQ$ of $G$ and $\hG$ 
containing respectively $B$ and $\hB$. We will explain how to find the well covering pairs 
in the case when $X=G/Q\times\hG/\hQ$.\\

\paragraphe
Let $\lambda$ be a one parameter subgroup of $T$ and so of $\hT$.
We can describe the fix point set $X^\lambda$:
$$
X^\lambda=\bigcup_{\begin{array}{c}
w\in W_{P(\lambda)}\backslash W/W_Q\\
\hw\in \hW_{\hP(\lambda)}\backslash\hW/\hW_\hQ
\end{array}
}
G^\lambda w Q/Q\times \hG^\lambda\hw\hQ/\hQ.
$$
For $(w,\hw)\in W_Q\backslash W/W_{P(\lambda)}\times
\hW_\hQ\backslash\hW/\hW_{\hP(\lambda)}$, we set 
$$
C(w,\hw)=G^\lambda w^{-1} Q/Q\times \hG^\lambda\hw^{-1}\hQ/\hQ.
$$
Note that, for later use, we have introduced $\ ^{-1}$ in this definition.
Note also that 
$$
C^+(w,\hw)= P(\lambda) w^{-1} Q/Q\times \hP(\lambda)\hw^{-1}\hQ/\hQ.
$$\\

\paragraphe
Let $P$ be a parabolic subgroup of $G$.
We consider the cohomology ring $\H^*(G/P,\ZZ)$ of $G/P$.
Here, we use singular cohomology with integers coefficients.
If $V$ is a closed subvariety of $G/P$, we denote by $[V]$ the Poincaré dual class 
in $\H^*(G/P,\ZZ)$.

\begin{nota}
The elements of $W/W_P$ will be denoted $w$ as elements of $W$:
in other word, the notation do not distinguish a class and a representative. 
Each times we use this abuse the reader has to check that the considered quantities does not
depend on the representative.
\end{nota}

If $w\in W/W_P$, we denote the corresponding Schubert variety 
$\overline{BwP/P}$ by $\Lambda_w^P$.
We denote by $[\Lambda_w^P]$ the Poincaré dual class in $\H^*(G/P,\ZZ)$ of $\Lambda_w$.
We also denote by $[{\rm pt}]$ the Poincaré dual class of the point; 
note that, $[{\rm pt}]=[\Lambda_e]$.
Let us recall that 
$$
\H^*(G/P,\ZZ)=\bigoplus_{w\in W/W_{P(\lambda)}}\ZZ[\Lambda_w^P].
$$ 
We use similar notation for $\hG/\hP$.

We now consider the case when $P=P(\lambda)$ and $\hP=\hP(\lambda)$.
Since $P=G\cap \hP(\lambda)$, $G/P(\lambda)$ identifies with the orbit by $G$ of
$\hP(\lambda)/\hP(\lambda)$ in $\hG/\hP(\lambda)$; 
let $\iota\::\:G/P(\lambda)\longto\hG/\hP(\lambda)$ denote this closed immersion.
The map $\iota$ induces a map $\iota^*$ in cohomology:
$$
\iota^*\::\:\H^*(\hG/\hP(\lambda),\ZZ)\longto\H^*(G/P(\lambda),\ZZ).
$$
  
To simplify notation, we set  $P=P(\lambda)$ and $\hP=\hP(\lambda)$.
\begin{lemma}
\label{lem:etabircohom}
For $(w,\hw)\in W_Q\backslash W/W_{P}\times
\hW_\hQ\backslash\hW/\hW_{\hP}$.

Then the following are equivalent:
\begin{enumerate}
\item the pair $(C(w,\hw),\lambda)$ is covering,
\item  $\iota^*([\overline{\hQ\hw\hP/\hP}]).
[\overline{QwP/P}]=[{\rm pt}]$.
\end{enumerate}

\end{lemma}

\begin{proof}
 Consider the map:
$$
\eta\::\:G\times_{P}C^+(w,\hw)\longto X.
$$
Since the characteristic of $\kk$ is zero, $\eta$ is birational if and only if for $x$ in an open subset of $X$, $\eta^{-1}(x)$ is reduced to a point. 
Consider the projection $\pi\::\:G\times_{P}C^+(w,\hw)\longto G/P$.
For any $x$ in $X$, $\pi$ induces an isomorphism from $\eta^{-1}(x)$ onto the following 
locally closed subvariety of $G/P$:
$F_x:=\{hP\in G/P\::\:h^{-1}x\in  C^+(w,\hw)\}$.

Let $(g,\hg)\in G\times\hG$ and set $x=(gQ/Q,\hg\hQ/\hQ)\in X$. 
We have:
$$
\begin{array}{rl}
F_x=&\{hP/P\in G/P\::\: h^{-1}gQ/Q\in Pw^{-1}Q/Q{\rm\ and\ }
h^{-1}\hg\hQ/\hQ\in \hP\hw^{-1}\hQ/\hQ\}\\
   =& \{hP/P\in G/P\::\: h^{-1}\in (Pw^{-1}Qg^{-1})\cap (\hP\hw^{-1}\hQ\hg^{-1})\}\\
   =& \iota(gQwP/P)\cap (\hg\hQ\hw\hP/\hP).
\end{array}
$$

Let us fix $g$ arbitrarily. 
By Kleiman's Theorem (see~\cite{Kleiman}), there exists an open subset of $\hg$'s in $\hG$  
such that the intersection
$ \overline{gQwP/P}\cap \overline{\hg\hQ\hw\hP/\hP}$ is transverse.
Moreover (see for example \cite{BK}), one may assume 
that $(gQwP/P)\cap (\hg\hQ\hw\hP/\hP)$ is dense in 
$\overline{gQwP/P}\cap \overline{\hg\hQ\hw\hP/\hP}$.
We deduce that the following are equivalent:
\begin{enumerate}
\item \label{cond1abo}for generic $\hg$, $F_x$ is reduced to a point,
\item $\iota^*([\overline{\hQ\hw\hP/\hP}]).[\overline{QwP/P}]=[{\rm pt}]$. 
\end{enumerate}

Since $\eta$ is $G$-equivariant, the above Condition~\ref{cond1abo} is clearly 
equivalent to the fact that $\eta$ is birational.
\end{proof}\\

\paragraphe
\begin{nota}
From now on, $\lg$ and $\lb$ will denote the Lie algebras of $G$ and $B$, 
$R$ will denote the set of roots of $\lg$ and $R^+$ those of  positive ones.
We denote by $\rho$ the half sum of the positive roots of $\lg$.
We will also use the following similar notation for $\hG$: $\hlg,\,\hlb,\,\hR,\,\hR^+,\,\hrho$.
\end{nota}

Let $w\in W/W_P$ and consider the associated $B$-orbit in $G/P$. 
We define $\gamma^P_w$ to be the sum of the weights of $T$ in the  normal space
at $wP/P$ of $BwP/P$ in $G/P$. 
Similarly, we define $\hat\gamma_\hw^\hP$.

\begin{lemma}
\label{lem:calculgamma}
We assume that $P$ is standard (that is, contain $B$, that is, $\lambda$ dominant).
We keep above notation. Then, $\gamma^P_w$ is the sum of the weights of $T$ in 
$\lg/(\lb+w.\lp)$.

Let $\tilde{w}$ be the longest element in the class $w\in W/W_P$.
We have:
$$
\gamma_w^{P}=-(\rho+ \tilde{w}\rho).
$$
\end{lemma}

\begin{proof}
Consider the map $G\longto G/P,\,g\mapsto gwP/P$ and its tangent map 
$\varphi\,:\,\lg\longto T_{wP/P}G/P$. Since $\varphi$ is $T$-equivariant, surjective and
$\lb+w.\lp=\varphi^{-1}(T_{wP/P}BwP/P)$,  $\gamma^P_w$ is the sum of the weights of $T$ in 
$\lg/(\lb+w.\lp)$.

Consider the $G$-equivariant surjective map $\pi\,:\,G/B\longto G/P$.
The definition of $\tilde{w}$ implies that $\pi(\tilde{w}B/B)=wP/P$ and 
$B\tilde{w}B/B$ is open in $\pi^{-1}(\overline{BwP/P})$.
It follows that $\gamma^P_w$ is  the sum of the weights of $T$ in the  normal space
at $\tilde{w}B/B$ of $B\tilde{w}B/B$ in $G/B$; that is,  the sum of the weights of $T$ in 
$\lg/(\lb+\tilde{w}.\lb)$. 
Since the sum of all the roots is zero, we obtain that
$$
\begin{array}{rl}
-\gamma^{P}_w&=\sum_{\alpha\in  R^+\cup \tilde{w}R^+}\alpha\\
&=\frac{1}{2}\left( \sum_{\alpha\in  R^+}\alpha+\sum_{\alpha\in  \tilde{w}R^+}\alpha+
\sum_{\alpha\in  R^+\backslash \tilde{w}R^+}\alpha+\sum_{\alpha\in \tilde{w}R^+\backslash R^+}\alpha
\right)\\
&=\frac{1}{2}\left(2\rho+2\tilde{w}\rho+\sum_{\alpha\in  R^+\backslash \tilde{w}R^+}\alpha+\sum_{\alpha\in \tilde{w}R^+\backslash R^+}\alpha
\right)
\\
&=\rho+\tilde{w}\rho
\end{array}
$$
where the last equality holds since $ R^+\backslash \tilde{w}R^+=-(\tilde{w}R^+\backslash R^+)$.
\end{proof}\\

\begin{remark}
In \cite{BK}, Belkale  and Kumar defined characters $\chi_{w^{-1}}$ for $w$ of minimal
length its coset in $W_P\backslash W$. 
We have 
$\gamma_{w^{-1}}^{B,P}=-\tilde{w}^{-1}w^P(\chi_{w^{-1}})$,
where $w^P$ denote the longest element of $W_P$.   
\end{remark}

\bigskip
\paragraphe
Let us fix again a dominant one parameter subgroup $\lambda$ of $T$, $w\in W$ and $\hw\in \hW$.

To simplify notation, we set $P=P(\lambda)$, $C=C(w,\hw)$ and $C^+=C^+(w,\hw)$.
Consider 
$$
\eta\::\:G\times_P C^+\longto X=G/Q\times\hG/\hQ.
$$

\begin{nota}
If $Y$ is a smooth variety of dimension $n$, $\Tau Y$ denotes its tangent bundle.
The line bundle $\bigwedge^n\Tau Y$ over $Y$ will be called the {\it determinant bundle} and
denoted by $\Det Y$.
If $\varphi\::\:Y\longto Y'$ is a morphism between smooth variety, we denote
by $T\varphi\::\:\Tau Y\longto\Tau Y'$ its tangent map, and by 
$\Det\varphi\::\:\Det Y\longto\Det Y'$ its determinant.\\
\end{nota}

Consider the restriction of $T\eta$ to $C^+$:
$$
T\eta_{|C^+}\::\:\Tau(G\times_PC^+)_{|C^+}
\longto \Tau(X)_{|C^+},
$$
and the restriction of $\Det\eta$ to $C^+$:
$$
\Det\eta_{|C^+}\::\:\Det(G\times_P C^+)_{|C^+}
\longto \Det(X)_{|C^+}.
$$
Since $\eta$ is $G$-equivariant, the morphism $\Det\eta_{|C^+}$ is $P$-equivariant;
it can be thought as a $P$-invariant section of the line bundle 
$\Det(G\times_P C^+)_{|C^+}^*\otimes\Det(X)_{|C^+}$ over $C^+$.
We denote by $\Li_{P,w,\hw}$ this last $P$-linearized line bundle on $C^+$.

\begin{lemma}
\label{lem:DetC}
We assume that $P,Q,\hP$ and $\hat{Q}$ are standard.
Let $(w,\hw)\in W/W_P\times\hW/\hW_\hP$ be such that 
$BwP/P$ and $\hB\hW\hP/\hP$ are open in $QwP/P$ and $\hQ\hW\hP/\hP$
respectively.
Let $\tilde{w}$ (resp. $\tilde{\hw}$) be the longest element in the class of $w$ (resp. $\hw$) 
in $W/W_P$ (resp. $\hW/\hW_\hP$).

Then, the torus $T$ acts on the fiber over the point $(\tilde{w}^{-1}Q/Q,\tilde{\hw}^{-1}\hQ/\hQ)$ 
in $\Li_{P,w,\hw}$ by the character
$$r_T(\tilde{\hw}^{-1}\gamma^{\hP}_\hw)+\tilde{w}^{-1}\gamma^{P}_w-\gamma^{P}_e.$$
\end{lemma}

\begin{proof}
If $Z$ is a locally closed subvariety of a variety $Y$ and $z$ is a point of $Z$,
we denote by $N_z^Y(Z)$ the quotient $T_zY/T_zZ$ of the tangent spaces at $z$ of $Y$ and $Z$.
If $V$ is a $T$-module $\St_T(V)$ denote the multiset of the weights of $T$ in $V$.
Let $\chi$ denote the character of the action of $T$ on the fiber over the point 
$x=(\tilde{w}^{-1}Q/Q,\tilde{\hw}^{-1}\hQ/\hQ)$ in $\Li_{P,w,\hw}$. 
Let $\lp$ denote the Lie algebra of $P$.

Since $\eta$ induces the identity on $C^+$ (canonically embedded in $G\times_PC^+$),
we have:
$$
\chi=-\sum_{\alpha\in \St_T(N_x^{G\times_PC^+}(C^+))}\alpha+
\sum_{\alpha\in \St_T(N_x^X(C^+))}\alpha.
$$
Moreover, we have the following $T$-equivariant isomorphisms:
$$\begin{array}{ll}
  N_x^{G\times_PC^+}(C^+)&\simeq \lg/\lp\simeq N_e^G(P)\simeq\lg/\lp,\\
N_x^X(C^+)&\simeq N_{\tilde{\hw}^{-1}\hQ/\hQ}^{\hG/\hQ}(\hP\hw^{-1}\hQ/\hQ)\oplus 
N_{\tilde{w}^{-1}Q/Q}^{G/Q}(Pw^{-1}Q/Q)\\
&\simeq N_{\tilde{\hw}^{-1}\hB/\hB}^{\hG/\hB}(\hP\tilde{\hw}^{-1}\hB/\hB)\oplus 
N_{\tilde{w}^{-1}B/B}^{G/B}(Pw^{-1}B/B)\\
&\simeq \hlg/(\hlp+\tilde{\hw}^{-1}\hlb)\oplus
\lg/(\lp+\tilde{w}^{-1}\lb)
\end{array}
$$
Now, the lemma is direct consequence of Lemma~\ref{lem:calculgamma}.
\end{proof}\\

\begin{nota}
 With notation of Lemma~\ref{lem:DetC}, we set $\theta_w^P=\tilde{w}^{-1}\gamma^{P}_w$
and $\theta^\hP_\hw=\tilde{\hw}^{-1}\gamma^{\hP}_\hw$.
\end{nota}

One can now describe the well covering pairs of $X$:

\begin{prop}
  \label{prop:wellcov}
Let $\lambda$ be a dominant one parameter subgroup of $T$.
Let $(w,\hw)\in W/W_P\times\hW/\hW_\hP$ be such that 
$BwP/P$ and $\hB\hat w\hP/\hP$ are open in $QwP/P$ and $\hQ\hat w\hP/\hP$
respectively.

The following are equivalent:
\begin{enumerate}
\item
The pair $(C(w,\hw),\lambda)$ is well covering.
\item \label{eq2:wellcov}
  \begin{tabular}{r@{}c@{}l@{}c@{}l@{}l}
$\iota^*([\Lambda_\hw^\hP])$&$.$&$[\Lambda_w^P]$&$-$&$[\Lambda_e^P]$&$=0$, and \\
$\langle \hw\lambda,\gamma^\hP_\hw\rangle$&$+$&$\langle w\lambda,\gamma^P_w\rangle$&$-$&
$\langle\lambda,\gamma^P_e\rangle$&$=0$.
  \end{tabular}
\end{enumerate}
\end{prop}

\begin{proof}
  By Lemma~\ref{lem:etabircohom}, we may (and shall) assume that
$(C(w,\hw),\lambda)$ is  covering.
Note that $w\lambda$ is well defined, since $W_P$ is precisely the stabilizer of $\lambda$ in $W$.
We chose $\tilde{w}$ and $\tilde{\hw}$ as in Lemma~\ref{lem:DetC}. Note that
$$\langle\lambda,r_T(\tilde{\hw}^{-1}\gamma^{\hP}_\hw)+\tilde{w}^{-1}\gamma^{P}_w-\gamma^{P}_e\rangle=
\langle \hw\lambda,\gamma^\hP_\hw\rangle+\langle w\lambda,\gamma^P_w\rangle-
\langle\lambda,\gamma^P_e\rangle.$$ 
In particular, by  Lemma~\ref{lem:DetC}, 
$\langle \hw\lambda,\gamma^\hP_\hw\rangle+\langle w\lambda,\gamma^P_w\rangle-
\langle\lambda,\gamma^P_e\rangle=0$  if and only if $\lambda$ acts trivially on the restriction of 
$\Li_{P,w,\hw}$ to $C(w,\hw)$.

Assume that $(C(w,\hw),\lambda)$ is well covering. Then $\Det\eta_{|C}$ is  non 
identically zero. Since $\Det\eta_{|C}$ is a $G^\lambda$-invariant section of ${\Li_{P,w,\hw}}_{|C}$,
$\lambda$ which fixes point wise $C$ has to act trivially on ${\Li_{P,w,\hw}}_{|C}$ (see for example
Lemma~\ref{lem:LGx}).

Conversely, assume Condition~\ref{eq2:wellcov} satisfied.
By Lemma~\ref{lem:DetC}, this implies that $\mu^{{\Li_{P,w,\hw}}_{|C}}(C,\lambda)=0$.
But, since $\eta$ is birational, $\Det\eta$ is $G$-invariant and non zero; hence, 
$\Det\eta_{|C^+}$ is $P$-invariant and non zero.
So, Proposition~\ref{prop:picC+} show that the restriction of $\Det\eta_{|C}$ is non identically zero.
Since $\eta$ is birational, this implies that $\eta$ is an isomorphism over an open subset 
intersecting $C$.
\end{proof}\\

\subsection{The case $X=G/B\times \hG/\hB$}
\label{sec:GBGB}

\paragraphe
We denote by $\lr(G,\hG)$ the cone of the pairs $(\nu,\hnu)\in \Chi(T)_\QQ\times\Chi(\hT)_\QQ$ 
such that for a positive integer $n$, $n\hnu$ and $n\nu$ are dominant weights such that 
$V_{n\nu}\otimes V_{n\hnu}$ contains non zero $G$-invariant vectors.

From now on, $X=G/B\times\hG/\hB$. 
By Proposition~\ref{prop:interp}, a point  $(\nu,\hnu)$ belongs to $\lr(G,\hG)$ if and only if
$r^{\Delta G}(\Li_{(\nu,\hnu)})$ belongs to $\tc^G(X)=\sac^G(X)$.\\

\paragraphe
Consider the $G$-module $\hlg/\lg$. 
Let $\chi_1,\cdots,\chi_n$ be the set of the non trivial weights of $T$ on $\hlg/\lg$.
For $I\subset\{1,\cdots,n\}$, we will denote by $T_I$ the neutral component of the intersection
of the kernels of the $\chi_i$'s with $i\in I$. 
A subtorus of the form $T_I$ is said to be {\it admissible}.
The subtorus $T_I$ is said to be {\it dominant} if $\Y(T_I)_\QQ$ is spanned by its
intersection with the dominant chamber of $Y(T)_\QQ$. 
Notice that the set of the $\chi_i$'s being stable by the action of $W$, any admissible subtorus is
conjugated by an element of $W$ to a dominant admissible subtorus. 
A one parameter subgroup of $T$ is said to be {\it admissible} if its image is.

To each $\chi_i$, we associate the hyperplane $\Hcal_i$ in $\Y(T)_\QQ$ spanned by the 
$\lambda\in\Y(T)$ such that $\chi\circ\lambda$ is trivial.
The $\Hcal_i$'s form an $W$-invariant arrangement of hyperplane in  $\Y(T)_\QQ$.
Moreover, $\Y(T_I)_\QQ$ is the intersection of the $\Hcal_i$'s with $i\in I$.\\

To simplify, in the following statement we assume that $\ac^G(X)$ has a non empty interior
in $\Pic^G(X)_\QQ$. In fact, this assumption is equivalent to say that no ideal of $\lg$ is an 
ideal of $\hlg$.

\begin{theo}
\label{th:ppal}
  We assume that no ideal of $\lg$ is an ideal of $\hlg$.
\begin{enumerate}
\item \label{ass:nonvide}
The interior of $\lr(G,\hG)$ in $\Chi(T)_\QQ\times\Chi(\hT)_\QQ$ is not empty.
\item \label{ass2:BeSj}
Let $\Face$ be a face of $\lr(G,\hG)$ of codimension $r$ which intersects the interior of the 
Weyl chamber. 
Then there exists a dominant  admissible subtorus $T_I$ (with $I\subset\{1,\cdots,n\}$) of $T$ 
of dimension $r$, a dominant indivisible one parameter subgroup $\lambda$ of $T_I$, 
and an irreducible component $C(w,\hw)$ of $X^\lambda$ (and $X^{T_I}$) such that:
\begin{enumerate}
\item \label{cond:H}
$\iota^*([\Lambda_\hw^{\hP(\lambda)}]).[\Lambda_w^{P(\lambda)}]=[\Lambda_e^{P(\lambda)}]
\in H^*(G/P(\lambda),\ZZ)$, 
\item \label{cond:theta}
$r_T(\htheta^{\hP(\lambda)}_\hw)+\theta^{P(\lambda)}_w-\theta^{P(\lambda)}_e$ is trivial on $T_I$, and
\item $\Face$ is the set of $(\nu,\hnu)\in\lr(G,\hG)$ such that 
$\langle w\lambda,\nu\rangle+\langle\hw\lambda,\hnu\rangle =0$.
\end{enumerate}

\item\label{ass3:BeSj}
Conversely, let $\lambda$ be a dominant one parameter subgroup of $T$ and $C=C(w,\hw)$ be an 
irreducible component $X^\lambda$. Set $I=\{i=1,\cdots,n\:|\:\chi\circ\lambda{\rm\ is\ trivial}\}$ 
and denote by $r$ the dimension of $T_I$.
If
\begin{enumerate}
\item $\iota^*([\Lambda^{\hP(\lambda)}_\hw]).[\Lambda^{P(\lambda)}_w]=[\Lambda^{P(\lambda)}_e]
\in H^*(G/P(\lambda),\ZZ)$ and 
\item $r_T(\htheta^{\hP(\lambda)}_\hw)+\theta^{P(\lambda)}_w-\theta^{P(\lambda)}_e$ is trivial on $T_I$,
\end{enumerate}
then the set of  $(\nu,\hnu)\in\lr(G,\hG)$ such that
$\langle w\lambda,\nu\rangle+\langle\hw\lambda,\hnu\rangle=0$  is a face of $\lr(G,\hG)$ 
of codimension $r$.
\end{enumerate}
\end{theo}

\begin{proof}
By \cite[Corollaire~1]{MontRess:LRcone}, the codimension of $\ac^G(X)$ in $\Pic^G(X)_\QQ$ is the dimension of the 
generic isotropy of $T$ acting on $\hlg/\lg$. This generic isotropy is also the kernel of the action of $T$ on $\hG/G$.
So, it is contained in $\bigcap_{\hg\in \hG}\hg G\hg^{-1}$. 
Since this group is distinguish in $\hG$ and $G$, it is finite by assumption.
Hence, the interior of $\ac^G(X)$ in $\Pic^G(X)_\QQ$ is non empty.\\

  Let $\Face$ be a face of $\ac^G(X)$ of codimension $r$.
Let $S$ be a torus in $G$, $\lambda$ be an indivisible  one parameter subgroup of $S$ 
and $C$ be an irreducible component of $X^S$ satisfying Corollary~\ref{cor:facecovGB}.
Up to conjugacy, we may assume that $S$ is contained in $T$ and $\lambda$ is dominant.

\noindent\underline{Claim 1.}
The subtorus $S$ is  admissible.

Notice that $C$ is isomorphic to the variety of the complete flags $G^S\times \hG^S$ 
and $S$ is the isotropy of a general point in $C$.
Hence, $S$ is the isotropy in $B^S$ of a  general point in $\hG^S/\hB^S$, 
that is, by Bruhat's Theorem, of a general point in $\hB^S/\hT^S$. 
Since $U^S$ acts freely on $\hB^S/\hT^S$, $S$ is conjugated to the isotropy in $T$ of a general 
point in $\hU^S/U^S$ which is isomorphic to $\hlu^S/\lu^S$. 
Finally, on may assume that $S$ is the isotropy of a general point in $\hlg^S/\lg^S=(\hlg/\lg)^S$; 
that is, since $S$ is diagonalizable, that $S$ is the kernel of the action of $T$ on $(\hlg/\lg)^S$. 
Then, $S$ is admissible.\\

Let $I\subset\{1,\cdots,n\}$ such that $S=T_I$.
Since $(C,\lambda)$ is well covering, Proposition~\ref{prop:wellcov} shows that 
$w,\,\hw$ and $\lambda$ satisfy Conditions~\ref{cond:H} and \ref{cond:theta} of 
the theorem. Assertion~\ref{ass2:BeSj} follows.\\

Conversely, let $w,\,\hw,\,\lambda$ and $I$ as in Assertion~\ref{ass3:BeSj} of the theorem.
By Proposition~\ref{prop:wellcov}, $(C(w,\hw),\lambda)$ is a well covering pair.
By Lemma~\ref{lem:muCcov}, the set of $\Li\in\tc^G(X)$ such that $\mu^\Li(C,\lambda)=0$ is a face  $\Face'$ of 
$\tc^G(X)$ eventually empty. It remains to prove that the codimension of $\Face'$ equal the dimension $r$ of $T_I$.

Consider the  linear map $\rho$ induced by the restriction:
$\rho\ :\ \Pic^G(X)_\QQ\longto\Pic^{G^\lambda}(C(w,\hw))_\QQ$.
We firstly prove the \\

\noindent\underline{Claim 3.} $\rho$ is surjective.

Let $x$ be a point in $C(w,\hw)$. Let $(\hB_1,B_1)$ denote the isotropy of $x$ in $G\times\hG$.
Let $\Li\in\Pic^{G^\lambda}(C(w,\hw))$. By Lemma~\ref{lem:PicGQ} applied to $C$, 
there exists characters $\hnu$ and $\nu$ of $\hB_1^\lambda$ and $B_1^\lambda$ and a positive integer $n$  
such that $\Li^\otimes$ is the restriction of the $G^\lambda\times\hG^\lambda$-linearized line bundle 
$\Li_{(\nu,\hnu)}$. But,  $\nu$ and $\hnu$ can be extended to characters $\nu'$ and $\hnu'$ of $B_1$ and $\hB_1$.
Then, $\Li^{\otimes n}$ is the image by $\rho$ of $r^{\Delta G}(\Li_{\nu',\hnu'})$.\\

By Claim~3 and Theorem~\ref{th:wc2face}, it remains to prove that 
$\tc^{G^\lambda}(C(w,\hw))$ has codimension $r$ in $\Pic^{G^\lambda}(C(w,\hw))_\QQ$.
By \cite[Corollaire~1]{MontRess:LRcone}, this codimension is the dimension of the kernel of the action of $T$ 
over $\hlg^\lambda/\lg^\lambda=(\hlg/\lg)^\lambda$. By definition $T_I$ is the neutral component of this kernel.
\end{proof}\\

\begin{coro}
\label{cor:face1}
  We assume that no ideal of $\lg$ is an ideal of $\hlg$.
\label{cor:BeSj}
Any dominant weight $(\nu,\hnu)$ belongs to  $\lr(G,\hG)$ if and only if
$$
\langle w\lambda,\nu\rangle+\langle \hw\lambda,\hnu\rangle \geq 0,
$$
for all indivisible dominant admissible one parameter  subgroup $\lambda$ of $T$ 
and for all $(w,\hw)\in W/W_{P(\lambda)}\times \hW/\hW_{\hP(\lambda)}$ such that
\begin{enumerate}
\item $\iota^*([\Lambda_\hw^{\hP(\lambda)}]).[\Lambda^{P(\lambda)}_w]=[\Lambda^{P(\lambda)}_e]
\in H^*(G/P(\lambda),\ZZ)$, and 
\item $\langle \hw\lambda,\gamma^{\hP(\lambda)}_\hw\rangle+
\langle w\lambda,\gamma^{P(\lambda)}_w\rangle=
\langle\lambda,\gamma^{P(\lambda)}_e\rangle$.
\end{enumerate}
Moreover, the above inequalities are pairwise distinct and no one 
can be omitted. 
\end{coro}

\begin{proof}
By  Proposition~\ref{prop:desc1}, the cone $\sac^G(X)$ is characterized as a part 
of the dominant chamber by the inequalities $\mu^\Li(C,\lambda)\leq 0$,  
for all well covering pair $(C,\lambda)$ with a dominant one parameter subgroup 
$\lambda$ of $T$.
Since $\sac^G(X)$ has non empty interior, it is sufficient to keep only the inequalities 
corresponding to faces of codimension one. 
By the theorem, all these inequalities are in the corollary.
The first part of the corollary is proved.

Consider an inequality $\langle w\lambda,\nu\rangle+\langle \hw\lambda,\hnu\rangle\geq 0$ as in the statement of the corollary.
By Theorem~\ref{th:ppal}, the set $(\nu,\hnu)$ such that 
$\langle w\lambda,\nu\rangle+\langle \hw\lambda,\hnu\rangle =0$ is a face $\Face$ 
of codimension one of $\lr(G,\hG)$. This inequality cannot be omitted unless $\Face$ is a face of
the dominant chamber in $\Chi(T)_\QQ\times\Chi(\hT)_\QQ$. It is easy to check that this is not possible.

It remains to prove that these inequalities are pairwise distinct.
But the stabilizer in $W$ (resp. $\hW$) of $\lambda\in Y(T)$ (resp.  $\lambda\in Y(\hT)$) is
precisely $W_{P(\lambda)}$ (resp. $\hW_{\hP(\lambda)}$).
It follows that the inequalities of the corollary are pairwise distinct.
\end{proof}

\paragraphe
Corollary~\ref{cor:face1} gives a minimal list of inequalities which determine $\lr(G,\hG)$ as a part of the
dominant chamber. In \cite{sjamaar}, Berenstein-Sjamaar gave another list containing redundant inaqualities.
An hyperplane given by one redundant inequality cannot intersect $\lr(G,\hG)$ along a face of codimension 
one containing  pairs of strictely dominant weights. 
We now want to understand better this intersection.

Let us fix an indivisible dominant admissible one parameter subgroup $\lambda$. 
Let  $(\Lambda_w,\,\Lambda_\hw)$ be a  pair of Schubert varieties  in $G/P(\lambda)$ and $\hG/\hP(\lambda)$ 
such that $\iota^*([\Lambda_\hw]).[\Lambda_w]=d.[{\rm pt}]\in H^*(G/P(\lambda),\ZZ)$ for some positive integer $d$.
By \cite{sjamaar}, for all $(\nu,\hnu)\in\lr(G,\hG)$, we have 
$\langle w\lambda, \nu \rangle+\langle\hw\lambda,\hnu\rangle\geq 0$.
In particular, the set $\Face$ of $(\nu,\hnu)\in\lr(G,\hG)$ such that
$\langle w\lambda, \nu \rangle+\langle\hw\lambda,\hnu\rangle=0$ is a face
of $\lr(G,\hG)$.

\begin{theo}
\label{th:redundantineq}
With above notation, we assume that $d\neq 1$ or that 
$\langle \hw\lambda,\gamma^{\hP(\lambda)}_\hw\rangle+
\langle w\lambda,\gamma^{P(\lambda)}_w\rangle\neq\langle\lambda,\gamma^{P(\lambda)}_e\rangle$.

Then, the face $\Face$ does not contain any weight $(\nu,\hnu)$ with $\nu$ strictly dominant.
\end{theo}

\begin{proof}
Set $x=(wB/B,\hw\hB/\hB)$.  
Let $C$ denote the irreducible component of $X^\lambda$ containing $x$.
By absurd, we assume that $\Face$ contains a weight $(\nu,\hnu)$ with $\nu$ strictly dominant.
According to Proposition~\ref{prop:wellcov}, it remains to prove that $(C,\lambda)$ is well
covering to obtain a contradiction.

Consider the parabolic subgroup of $\hG$ containing $\hB$ such that $\Li_{(\nu,\hnu)}$ is an
ample line bundle on $G/B\times \hG/\hQ$ (denoted $\overline{X}$ from now on). 
Consider the natural $G\times\hG$-equivariant morphism $p\,:\,X\longto \overline{X}$.
Set $\overline{x}=p(x)$. 
Let $\overline{C}$ denote the  irreducible component of $\overline{X}^\lambda$ containing $\overline{x}$
and $\overline{C}^+$ the corresponding Bialinicki-Birula cell. 

Since $\overline{C}$ (resp. $\overline{C}^+$) is an orbit of $G^\lambda\times\hG^\lambda$
(resp. $P(\lambda)\times\hP(\lambda)$), we have $p(C)=\overline{C}$ and 
$p(C^+)=\overline{C}^+$.
Since $\iota^*([\Lambda_\hw]).[\Lambda_w]=d.[{\rm pt}]$, the proof of Lemma~\ref{lem:etabircohom} 
shows that $(C,\lambda)$ is dominant. 
We deduce that $(\overline{C},\lambda)$ is dominant.
Now, Theorem~\ref{th:CxcGB} implies that  $(\overline{C},\lambda)$ is well covering.

Let $y\in C$ general and $g\in G$ such that $g^{-1}y\in C^+$.
Since $p(C)=\overline{C}$, $p(y)$ is general in $\overline{C}$.
But, $(\overline{C},\lambda)$ is well covering and $g^{-1}p(y)\in\overline{C}^+$; so,
$g\in P(\lambda)$. This proves that $(C,\lambda)$ is well covering.
\end{proof}

\subsection{Application to the tensor product}

\paragraphe
In this section, $G$ is assumed to be semisimple. 
We also fix an integer $s\geq 2$ and set $\hG=G^s$, $\hT=T^s$ and $\hB=B^s$.
We embed $G$ diagonally in $\hG$. 
Then $\sac^G(X)\cap \Chi(T)^{s+1}$ identifies with the $(s+1)$-uple 
$(\nu_1,\cdots,\nu_{s+1})\in\Chi(T)^{s+1}$ such that the for $n$ big enough $n\nu_i$'s are 
dominant weights and $V_{n\nu_1}\otimes\cdots\otimes V_{n\nu_{s+1}}$ contains a non zero 
$G$-invariant vector.

The set of weights of $T$ in $\hlg/\lg$ is simply the root system $\Phi$ of $G$. 
Let $\Delta$ be the set of simple roots of $G$ for $T\subset B$.
Let $I$ be a part $I$ of $\Delta$. Let $L(I)$ denote the Levi subgroup of $G$ containing $T$ and
having $\Delta-I$ as its simple roots. 
Let $T_I$ denote the neutral component of the center of $L(I)$; $T_I$  is  dominant.
Note that the dimension of $T_I$ is the cardinality of $I$. 
Moreover, all dominant admissible subtorus of $T$ is obtained in such a way.
We will also denote by $P(I)$ the standard parabolic subgroup with Levi subgroup $L(I)$.
We denote by $W_I$ the Weyl group $L(I)$.

Let $\lambda$ be a dominant one parameter subgroup of $T$.
For $(w_1,\hw)=(w_1,\cdots,w_{s+1})\in W\times \hW=W^{s+1}$, 
and $(\nu,\hnu)=(\nu_1,\cdots,\nu_{s+1})\in\Pic^G(X)_\QQ=\Chi(T)^{s+1}_\QQ$ we have:
\begin{itemize}
\item $r_T(\hw\hnu)=\sum_{i=2}^{s+1} w_i\nu_i$, and $r_T(\htheta^{\hP(\lambda)}_\hw)=
\sum_{i=2}^{s+1} \theta^{P(\lambda)}_{w_i}$,
\item $\iota^*([\Lambda^{\hP(\lambda)}_\hw])= [\Lambda^{P(\lambda)}_{w_2}]{\bf\cdot}\,\cdots\, {\bf\cdot}[\Lambda^{P(\lambda)}_{w_{s+1}}]$,
\end{itemize}

In \cite{BK}, Belkale and Kumar defined a new product denoted $\kbprod$ on the cohomology
groups $H^*(G/P,\ZZ)$ for any parabolic subgroup $P$ of $G$.
By Proposition~17 of \cite{BK}, this product $\kbprod$ has the following very interesting property.

For $w_i\in W/W_{P(I)}$, the following are equivalent:
\begin{enumerate}
\item $[\Lambda^{P(I)}_{w_1}].\cdots.[\Lambda^{P(I)}_{w_{s+1}}]=[\Lambda^{P(I)}_e]$ and,\\
the restriction of $\theta^{P(I)}_{w_1}+\cdots +\theta^{P(I)}_{w_{s+1}}-\theta^{P(I)}_e$ to $T_I$ is trivial;
\item $[\Lambda^{P(I)}_{w_1}]\kbprod\cdots\kbprod[\Lambda^{P(I)}_{w_{s+1}}]=[\Lambda^{P(I)}_e]$.
\end{enumerate}

Using this result of Belkale and Kumar our Theorem~\ref{th:ppal} gives the following corollary.
If $\alpha$ is a root of $G$, $\alpha^\vee$ denote the corresponding coroot.
If $\alpha$ is a simple root, $\omega_{\alpha^\vee}$ denote the corresponding fundamental weight.
  
\begin{coro}
\label{th:tensor}
\begin{enumerate}
\item \label{ass1:tensor}
A point $(\nu_1,\cdots,\nu_{s+1})\in\Chi(T)_\QQ^{s+1}$ belongs to the cone 
$\lr(G,G^s)$ 
if and only if 
\begin{enumerate}
\item \label{eq:dom}
each $\nu_i$ is dominant; that is $\langle\alpha^\vee,\nu_i \rangle \geq 0$ for all simple root $\alpha$.
\item \label{eq:mu}
for all simple root $\alpha$;
for all $(w_1,\cdots,w_{s+1})\in (W/W_{P(\alpha)})^{s+1}$ such that 
$[\Lambda^{P(\alpha)}_{w_1}]\kbprod\cdots\kbprod[\Lambda^{P(\alpha)}_{w_{s+1}}]=[\Lambda^{P(\alpha)}_e]
\in H^*(G/P(\alpha),\ZZ)$, we have:
$$
\sum_i\langle w_i\omega_{\alpha^\vee},\nu_i \rangle \geq 0.
$$
\end{enumerate}
\item We assume either that $s\geq 3$ or $\lg$ does not contain any factor of rank one.
In the above description of $\lr(G,G^s)$, the inequalities are pairwise distinct and no one can be omitted 
(neither in \ref{eq:dom} nor \ref{eq:mu}). 

\item \label{ass2:tensor}
Let $\Face$ be a  face of $\lr(G,G^s)$ of codimension $d$ which intersects the interior of
the dominant chamber.
There exist a subset $I$ of $d$ simple roots and
and $(w_1,\cdots,w_{s+1})\in (W/W_I)^{s+1}$  such that:
\begin{enumerate}
\item $[\Lambda^{P(I)}_{w_1}]\kbprod\cdots\kbprod[\Lambda^{P(I)}_{w_{s+1}}]=[\Lambda^{P(I)}_e]
\in H^*(G/P(I),\ZZ)$,
\item the subspace spanned by $\Face$ is the set 
$(\nu_1,\cdots,\nu_{s+1})\in\Chi(T)_\QQ^{s+1}$ such that:
$$
\forall \alpha\in I\hspace{3ex}\sum_i\langle w_i\omega_{\alpha^\vee},\nu_i \rangle =0.
$$
\end{enumerate}
\item \label{ass3:tensor}
Conversely, let $I$ be a subset of $d$ simple roots
and $(w_1,\cdots,w_{s+1})\in (W/W_I)^{s+1}$  such that
$[\Lambda^{P(I)}_{w_1}]\kbprod\cdots\kbprod[\Lambda^{P(I)}_{w_{s+1}}]=[\Lambda^{P(I)}_e]
\in H^*(G/P(I),\ZZ)$.
Then, the set of $(\nu_1,\cdots,\nu_{s+1})\in \lr(G,G^s)$ such that 
$$
\forall \alpha\in I\hspace{3ex}\sum_i\langle w_i\omega_{\alpha^\vee},\nu_i \rangle =0,
$$
is a face of codimension $d$ of $\lr(G,G^s)$.
\end{enumerate}
\end{coro}

\begin{proof}
Equations~\ref{eq:dom} are all different and are not repeated in Equations~\ref{eq:mu}.
Moreover, by \cite[Proposition~7]{MontRess:LRcone} they define codimension one faces of 
$\lr(G,G^s)$.

The rest of the corollary is a simple rephrasing of Theorem~\ref{th:ppal} and 
Corollary~\ref{cor:BeSj}.
\end{proof}\\

\begin{remark}
The description of the smaller faces of $\cone^G((G/B)^{s+1})$ gives an application of the 
Belkale-Kumar product $\kbprod$ for all the complete homogeneous spaces.  
\end{remark}

\bibliographystyle{amsalpha}
\bibliography{biblio}

\begin{center}
  -\hspace{1em}$\diamondsuit$\hspace{1em}-
\end{center}

\vspace{5mm}
\begin{flushleft}
N. R.\\
Universit{\'e} Montpellier II\\
D{\'e}partement de Math{\'e}matiques\\
Case courrier 051-Place Eug{\`e}ne Bataillon\\
34095 Montpellier Cedex 5\\
France\\
e-mail:~{\tt ressayre@math.univ-montp2.fr}  
\end{flushleft}

\end{document}